\documentclass[12pt]{article}

\usepackage{multirow,array}
\usepackage{graphicx,color,rotating}
\usepackage{amssymb, amsmath, amsfonts, bm}
\usepackage{psfrag}
\usepackage{epstopdf}
\usepackage{subfigure}
\usepackage{setspace}
\DeclareGraphicsExtensions{.pdf,.png,.gif,.jpg,.eps,.ps}

\def\I{{\mathcal{I}}}
\def\J{{\mathcal{J}}}
\def\C{{\mathcal{C}}}

\def\O{{\mathcal{O}}}

\def\X{{\mathcal{X}}}

\def\RR{{\mathbb{R}}}

\def\bA{{\bf A}}

\def\bF{{\bf F}}

\def\bI{{\bf I}}
\def\bJ{{\bf J}}

\def\bL{{\bf L}}

\def\ta{{\tilde{a}}}
\def\tr{{\tilde{r}}}

\def\a{{\mathbf a}}

\def\x{{\bm x}}

\def\u{{\bm u}}

\def\f{{\mathbf f}}

\def\0{{\mathbf 0}}
\def\bxi{\boldsymbol{\xi}}
\def\bnabla{\boldsymbol{\nabla}}

\def\Dpartial#1#2{ {\partial #1 \over \partial #2} }
\def\Dpartialmix#1#2#3{ {\partial^2 #1 \over \partial #2\,\partial #3} }

\def\Bmp#1{ \begin{minipage}{#1} }
\def\Emp{ \end{minipage} }
\def\Bmpc#1{ \begin{minipage}[c]{#1} }
\def\Bmpt#1{ \begin{minipage}[t]{#1} }
\def\Bmpb#1{ \begin{minipage}[b]{#1} }

\newcommand{\argmin}{\operatorname{argmin}}

\newtheorem{assumption}{Assumption}

\makeatletter

\newcommand{\Rmnum}[1]{\expandafter\@slowromancap\romannumeral #1@}
\makeatother

\title{An Optimal Model Identification For Oscillatory Dynamics With a
  Stable Limit Cycle}
\author{
Bartosz Protas\thanks{Department of Mathematics and Statistics, McMaster
  University, Hamilton, ON, Canada} \and 
Bernd R.~Noack\thanks{Institut PPRIME,
CNRS -- Universit\'e de Poitiers -- ENSMA, UPR 3346,
D\'epartement Fluides, Thermique, Combustion,
CEAT, 43 rue de l'A\'erodrome,
F-86036 POITIERS Cedex, France} \and 
Marek Morzy\'nski\thanks{Pozna\'n University of Technology,
Institute of Combustion Engines and Transportation, 
ul.\ Piotrowo 3, PL 60-965 Pozna\'n, Poland}}

\begin{document}
\maketitle

%\newcommand{\slugmaster}{%
%\slugger{siads}{xxxx}{xx}{x}{x--x}}%slugger should be set to juq, siads, sifin, or siims

%----abstract-------------------------------------------------------------------------------
\begin{abstract}
  We propose a general framework for parameter-free identification of
  a class of dynamical systems.  Here, the propagator is approximated
  in terms of an arbitrary function of the state, in contrast to a
  polynomial or Galerkin expansion {used in traditional
    approaches}.  The {proposed formulation} relies on
  variational data assimilation using measurement data combined with
  assumptions on the smoothness of the propagator.  This approach is
  {illustrated using} a generalized dynamic model describing
  oscillatory transients from an unstable fixed point to {a}
  stable limit cycle and arising in nonlinear stability analysis
  {as an example. This} 3-state model comprises an evolution
  equation for the dominant oscillation and an algebraic manifold for
  the low- and high-frequency components in an autonomous descriptor
  system.  The proposed optimal model identification technique employs
  mode amplitudes of the transient vortex shedding in a cylinder wake
  {flow} as {example measurements}.  The
  {reconstruction obtained with our technique features} distinct
  {and} systematic improvements over the well-known mean-field
  (Landau) model of the Hopf bifurcation.  The computational aspect of
  the identification method is thoroughly validated showing that good
  reconstructions can also be obtained in the absence of of accurate
  initial approximations.
\end{abstract}

%\begin{keywords} 

%\begin{AMS}
% 93A30   	Mathematical modeling (models of systems, model-matching, etc.)
% 65K10   	Optimization and variational techniques [See also 49Mxx, 93B40]
% 76D25   	Wakes and jets
%\end{AMS}

\pagestyle{myheadings}
\thispagestyle{plain}
%\markboth{B.~Protas, B.~R.~Noack And M.~Morzy\'nski}{Optimal Model Identification}

%----ToC-------------------------------------------------------------------------------
%\tableofcontents

%---- Main part --------------------------------------------------------
%\input{s1.tex}
%\input{s2.tex}
%\input{s3.tex} 
%\input{s4.tex}
%\input{s5A.tex}
%\input{s5B.tex}
%\input{s6.tex}%***********************************************************************
\section{Introduction}
\label{ToC:Introduction}
\begin{figure}
\begin{center}
\mbox{
\subfigure[]{\includegraphics[width=0.35\textwidth]{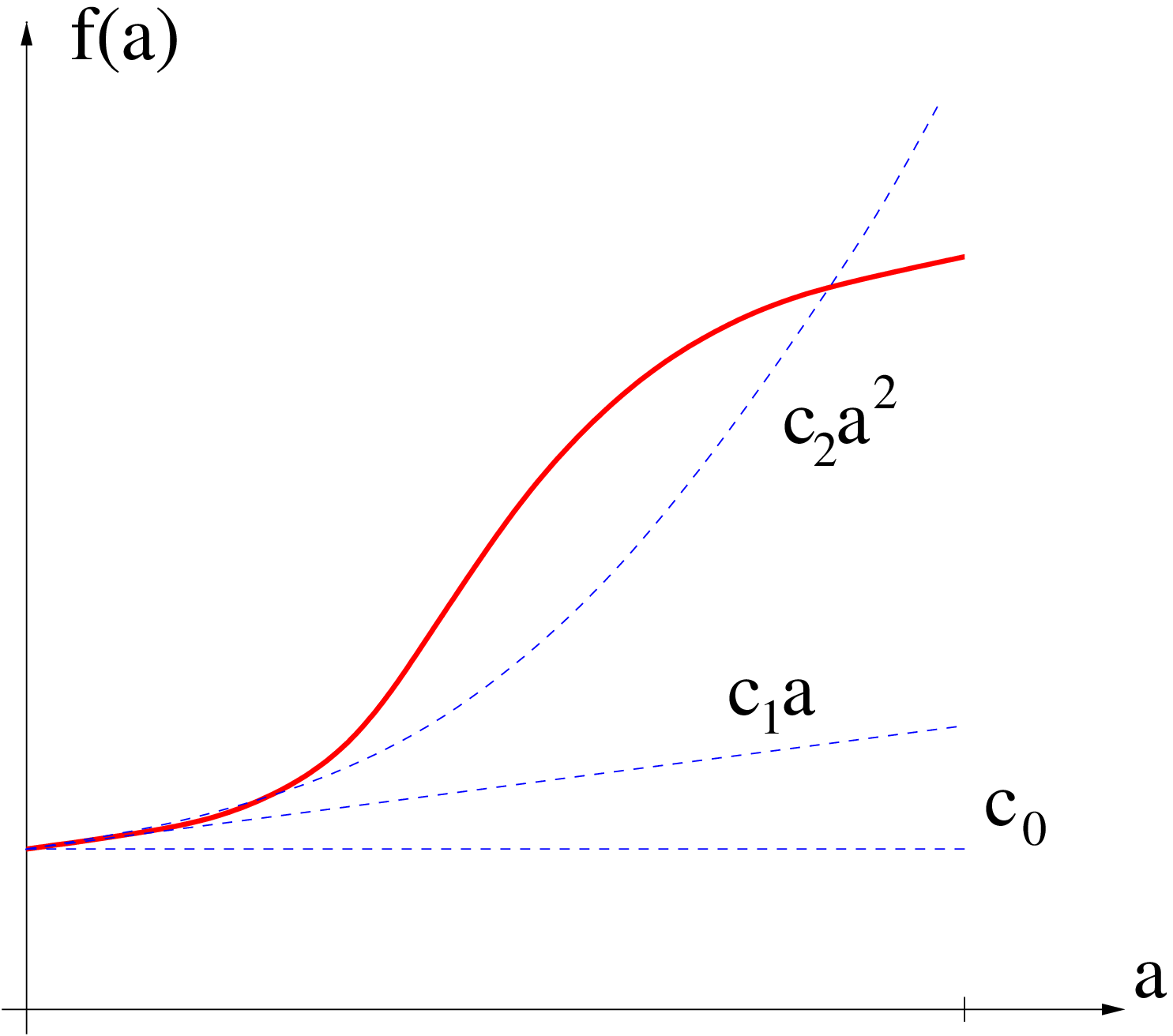}}
%\hfill
\qquad\qquad\qquad\qquad
\subfigure[]{\includegraphics[width=0.35\textwidth]{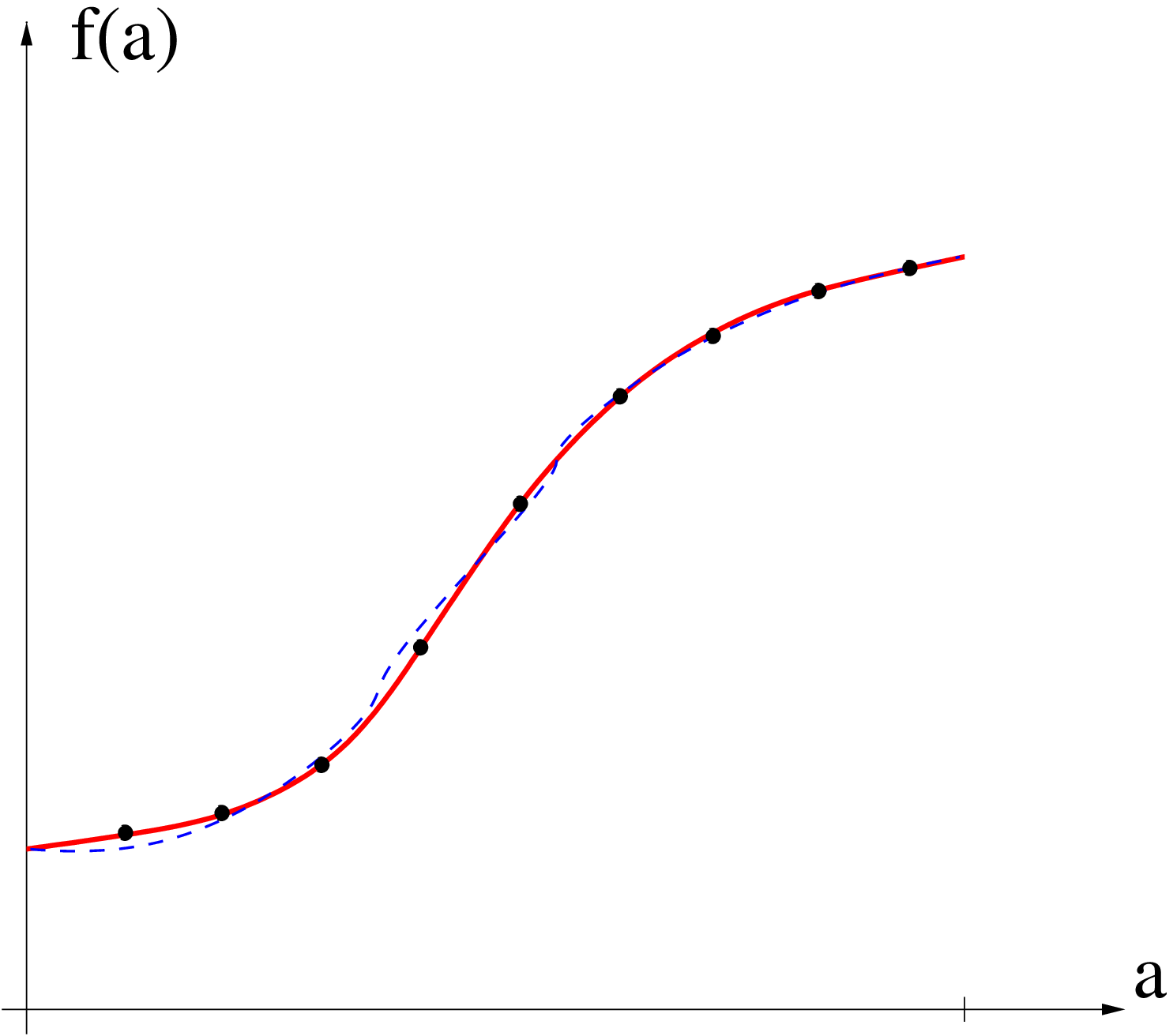}}}
{\mbox{
\Bmpt{0.4525\textwidth}
\centering Standard Approach:
\begin{itemize}
\item Step 1 --- select model dimension $N$,
\item Step 2 --- identify structure, e.g., 
\begin{equation*}
\frac{d}{dt}a \approx \sum_{k=0}^N c_k \, a^k, \hfill \hspace*{1.5cm} (\star)
\end{equation*}
\item Step 3 --- identify parameters $\{c_k\}_{k=0}^N$.
\end{itemize}
\Emp
\qquad
\Bmpt{0.4525\textwidth}
\centering Proposed Approach:
\begin{itemize}
\item Step 1 --- select model dimension $N$,
\item Step 2 --- perform parameter-free reconstruction of $f(a)$
  using the variational technique described below,
\item Step 3 --- truncate the reconstructed function to a suitable
  form such as, e.g., ($\star$).
\end{itemize}
\Emp
}}
\caption{{Schematic diagram illustrating the main idea behind
    (a) the standard approach and (b) the proposed new approach to
    model identification. In the figures the solid line represents the
    ``true'' descriptor function $f(a)$, whereas the dashed lines
    denote its reconstructions obtained with the two approaches. In
    Figure (b) the solid symbols represent the grid points used for
    the discretization of the variational formulation.}}
\label{fig:f}
\end{center}
\end{figure}

In this study we consider the problem of parameter-free identification
of a class of dynamical systems. The approach we propose is derived
from a general method for the reconstruction of the constitutive
relations in systems described by partial differential equations
(PDEs) which was initially introduced in \cite{bvp10} and further
developed in \cite{bp11a}, see also \cite{b12}. The idea is that,
given an autonomous evolution equation $\frac{d}{dt}\bm{a} = \bm{f} (
\bm{a} )$ for some quantity $\a(t)$ (defined in a finite or infinite
dimension with $t\ge 0$ denoting time), one seeks to optimally
reconstruct the flux function $\f(\a)$, so that the system outputs
best match, in a suitably defined sense, {with} the measurements
available.  {As shown schematically in Figure \ref{fig:f},} the
originality of this approach is that the function $\f(\a)$ is
reconstructed directly as a {\em continuous} object, rather than
employing a truncated polynomial or other Galerkin expansion.  In
\cite{bvp10,bp11a,b12} we reviewed the mathematical foundations and
some computational aspects of this approach applied to the
reconstruction of state-dependent transport coefficients in a class of
systems described by PDEs. In the present investigation we adapt this
method to the problem of identification of the propagator function
$\bm{f} \; : \; \RR^N \rightarrow \RR^N$, $\bm{f} =
[f_1,\ldots,f_N]^T$, in a finite-dimensional dynamical system
$\frac{d}{dt}\bm{a} = \bm{f} ( \bm{a} )$ with the state vector
$\bm{a}(t) = [a_1(t), \ldots, a_N(t)]^T \in \RR^N$. To fix attention,
instead of working with an abstract formulation, we will focus on a
specific problem concerning a three-dimensional ($N=3$) model
characterized by oscillatory dynamics with a stable limit cycle. In
many real-life physical problems with infinite-dimensional state
spaces and evolution described by PDEs such systems arise as
reduced-order models providing low-dimensional approximate description
in the neighborhood of fixed points and limit cycles.

In some cases propagator $\f$ may be derived from the description of
the full plant based on first principles.  In general, however, the
reduced-order model is inferred from experimental or numerical data.
Such reduced models are of paramount importance as test-beds for the
development of our understanding of system dynamics. They are also
useful as low-cost surrogates guiding optimization and real-time
control design for expensive full-scale models. In all cases, typical
model identification is generally performed in three steps: (1)
selection of the state space $\RR^N$ which is large enough to capture
the behaviour of interest and at the same time sufficiently small to
allow one to exploit the analytical/numerical advantages of the
surrogate plant; (2) structure identification of the propagator
$\bm{f}$, e.g., determination of the polynomial degree of its
components $f_i$, $i=1,\dots,N$; and (3) parameter identification,
e.g., inference of the polynomial coefficients from time-resolved
trajectories $t \mapsto \bm{a}(t)$. {These different steps are
  illustrated schematically in Figure \ref{fig:f}a.}  The challenge of
this approach is to find the right balance between the robustness of
the model identification, requiring only a few tunable parameters, and
a good accuracy for which a larger number of parameters is typically
needed.  Identification problems {can also be} solved using
variational techniques and this is the approach we will pursue in the
present study, {cf.~Figure \ref{fig:f}b}. Similar methods have
been developed for a broad range of problems in both the finite and
infinite-dimensional setting, including flow control in fluid
mechanics \cite{g03}, data assimilation in meteorology \cite{n97,k03}
and geophysics \cite{t05} to mention just a few application areas. The
related problem of state estimation is usually solved using various
filtering approaches such as the Kalman filter \cite{s94}.

We introduce now our model. The oscillatory fluctuation is
parameterized by
\begin{equation}
a_1 + \imath a_2 = r \> \exp(\imath \theta),
\label{eq:aexp}
\end{equation}
where $\imath$ is the imaginary unit, $r :=
\sqrt{a_1^2+a_2^2}$ is the amplitude of the fluctuation (``$:=$'' means ``equal to
by definition'') and $\theta
:= \arctan(a_2/a_1)$ the corresponding phase, 
while the base-flow deformation is characterized by a single parameter $a_3$.
Following {the} mean-field theory \cite{Stuart1971arfm},
we make the following assumptions about the structure of the system:
\begin{assumption}
\Bmp{0.1cm}\Emp 
\begin{enumerate}
% BRN: There is an unusal '.' after (i).
% BRN: Moreover, the more ususal numbering order is I->A->1->a->i
% BRN: Hence, I would expect a), b) instead of i) and ii).
\renewcommand{\theenumi}{(\alph{enumi})} 
\item In the plane $(a_1,a_2)$ the system exhibits an unstable fixed
  point at the origin and an attracting limit cycle,
\item the state variable $a_3$ is ``slaved'' to $a_1$ and $a_2$, i.e.,
  $a_3 = a_3(a_1,a_2)$, and
\item the dynamics is phase-invariant, 
  i.e., the descriptor system depends
  only on $r$.
% BRN: The phase invariance has to be the last assumption
% BRN: since it effects the dynamical system AND the manifold.
% vector field -> descriptor system to include the manifold
\end{enumerate}
\label{ass1}
\end{assumption}

\smallskip As regards the time $t$, we will assume that $t \in [0,T]$
for some $T>0$. As a general form of a dynamical system consistent
with Assumption \ref{ass1} we will consider
\begin{subequations}
\label{eq:DescriptorSystemB}
\begin{align}
\dot r{(t)}      &= g_1(r{(t)}) \> r{(t)},
\label{eq:DescriptorSystemB_1}
\\
\dot \theta{(t)} &= g_2(r{(t)}) ,
\label{eq:DescriptorSystemB_2}
\\
a_3{(t)}  &= g_3(r{(t)}),
\label{eq:DescriptorSystemB_3}
\end{align}
\end{subequations}
or, equivalently
\begin{subequations}
\label{eq:DescriptorSystemA}
\begin{align}
\frac{d}{dt} 
\begin{bmatrix} a_1(t) \\ a_2(t) \end{bmatrix}
& = \left( g_1(r) \, \bI + g_2(r) \, \bJ \right) 
\begin{bmatrix} a_1(t) \\ a_2(t) \end{bmatrix}
=: \begin{bmatrix} f_1 \\ f_2 \end{bmatrix},
\label{eq:DescriptorSystemA_1} \\
a_3(t) & = g_3(r(t)),
\label{eq:DescriptorSystemA_2}
\end{align}
\end{subequations}
where $\bI = \begin{bmatrix} 1 & 0 \\ 0 & 1 \end{bmatrix}$ and $\bJ =
\begin{bmatrix} 0 & -1 \\ 1 & 0 \end{bmatrix}$. Subsequently, we will
also use the notation $\bxi := [a_1, \ a_2]^T$ and $\I := [0,\
r_{max}]$, where $r_{max} := \sup_{t \in [0,T]} r(t)$.  {We note
  that the equations governing $a_1$ and $a_2$, or equivalently $r$
  and $\theta$, do not depend on $a_3$.} {Equations
  \eqref{eq:DescriptorSystemB_3} and \eqref{eq:DescriptorSystemA_2}
  describe the dependency of a slow variable on the fluctuation
  amplitude, and the usefulness of these algebraic equations will
  become clear in the context of the mean-field models discussed in
  Section \ref{ToC:MeanFieldSystem}.}  In \eqref{eq:DescriptorSystemB}
and \eqref{eq:DescriptorSystemA} the functions $g_i \; : \; \I
\rightarrow \RR$, $i=1,2,3$, are assumed sufficiently regular to make
the systems well posed (the question of the regularity of functions
$g_i$, $i=1,2,3$, will play an important role in our approach and will
be addressed in detail further below).  Without loss of generality, we
will also assume that $\dot\theta \ge 0$. System
\eqref{eq:DescriptorSystemA_1}, or
\eqref{eq:DescriptorSystemB_1}--\eqref{eq:DescriptorSystemB_2}, is
complemented with the initial condition, respectively, $\bxi(0) =
\bxi^0 := [a_1^0, \ a_2^0]^T$ and $r(0) = \Vert \bxi^0 \Vert$,
$\theta(0) = \arctan(a_2^0 / a_1^0)$. Dynamical systems of the type
\eqref{eq:DescriptorSystemB} or \eqref{eq:DescriptorSystemA} arise
commonly as a result of various rigorous and empirical model-reduction
strategies in diverse application areas such as fluid mechanics
\cite{dlgf94}, thermodynamics \cite{ll80} and phase transitions
\cite{y87}. The main contribution of this work is development of a
computational technique allowing one to reconstruct functions $g_1$,
$g_2$ and $g_3$ in \eqref{eq:DescriptorSystemB} in a very general form
based on some measurements. The key idea is to formulate a
least-squares minimization problem in which one of the functions
$g_i$, $i=1,2,3$, is the control variable. Then, a variational
gradient-based approach can be employed to find the optimal solution
in a suitable function space.

In the present investigation, we will consider system
\eqref{eq:DescriptorSystemB} as a reduced-order model of hydrodynamic
instabilities in open shear flows obtained using a suitable Galerkin
projection --- to fix attention without losing generality.  Although
this problem is very well-studied, our method provides a systematic
refinement of the state-of-the-art mean-field model, which is the
second main contribution of this study. Details of this problem are
introduced in the next Section. In the following Section we develop
our model identification approach, whereas in Section
\ref{ToC:Results} we present a number of computational results
concerning identification of the model for the system considered in
Section \ref{ToC:ModelIdentification}. Then, in Section
\ref{ToC:Discussion}, we analyze the computational performance of our
method and discuss the improvements it offers over the predictions of
some standard models applied to the problem in question. Summary and
conclusions are deferred to Section \ref{ToC:Conclusions}, whereas in
Appendices \ref{ToC:POD}, \ref{ToC:Regularity} and
\ref{ToC:Perturbation} we collect some technical results.

%***********************************************************************
%\section{Model identification problem /BERND}
\section{Example Problem --- Model Identification For a Vortex Shedding Instability}
\label{ToC:ModelIdentification}

In this Section we define a model identification problem associated
with the transient two-dimensional (2D) cylinder wake which is a
well-studied hydrodynamic instability
\cite{Jackson1987jfm,Morzynski1999cmame,Noack1994jfm}.  This flow is a
representative example of phenomena characterized by the Hopf
bifurcation with an unstable fixed point and a stable limit cycle.
{Additional} examples {in} this category {include} the Rossiter modes
of the flow over a cavity \cite{Rowley2006arfm} and other shear flows
\cite{Schmid2001book}.  First, in Section \ref{ToC:Wake}, we describe
the initial-boundary value problem for the infinite-dimensional
Navier-Stokes equation which is a system of coupled PDEs representing
the conservation of mass and momentum in the motion of viscous
incompressible fluid.  Then, in Section \ref{ToC:GalerkinExpansion}, a
low-dimensional Galerkin expansion is recalled which reduces the
kinematic description down to three modes whose amplitudes serve as
the state variables $a_1$, $a_2$ and $a_3$.  It will be demonstrated
that the system governing these variables satisfies in fact
Assumptions \ref{ass1} and is in the form
\eqref{eq:DescriptorSystemA}. Special forms of this system arising {as
  models} in numerous applications are discussed in Section
\ref{ToC:MeanFieldSystem}, whereas in Section \ref{ToC:Measurements}
we describe the measurements used as the basis for the
reconstructions.

%=======================================================================
\subsection{Cylinder Wake Flow}
\label{ToC:Wake}
The 2D flow around a circular cylinder is described here in the
Cartesian coordinate system $\bm{x}:=(x,y)$.  The origin $\bm{0}$
coincides with the center of the cylinder, the $x$-coordinate points
in streamwise direction, while $y$ represents the transverse
coordinate.  The velocity field $\bm{u}:=(u,v)$ has components $u$ and
$v$ aligned with the $x$-axis and $y$-axis, respectively, whereas $p$
represents the pressure field. The oncoming flow velocity is denoted
by $U$ and the cylinder diameter by $D$.  The Newtonian fluid is
characterized by uniform density $\rho$ and kinematic viscosity $\nu$.
The flow properties depend on the Reynolds number $Re := U \>D / \nu$.
Here, we consider the case of $Re=100$ which is far above the critical
Reynolds number of $47$ characterizing the onset of vortex shedding
(i.e., the Hopf bifurcation) \cite{Jackson1987jfm,Zebib1987jem}, and
far below the transitional Reynolds number of 187 which marks the
onset of three-dimensional instabilities \cite{Zhang1995pf}.

In the following, all quantities are assumed to be
non-dimensionalized with $U$, $D$   and $\rho$.  
The flow is considered in a rectangular domain surrounding the cylinder
\begin{equation}
\Omega := \{ (x,y) \> \colon \> x^2+y^2 \geq 1/4 \> 
\wedge \> -5 \leq x \leq 15 \> \wedge \> \vert y \vert \leq 5 \}
\label{eq:Omega}
\end{equation}
{and its} evolution is described by the incompressible Navier-Stokes equation
\begin{subequations}
\label{Eqn:IncompressibleFlow}
\begin{alignat}{2}
\label{Eqn:MassBalance}
\nabla \cdot \bm{u} &= 0&  & \textrm{in} \ \Omega,
\\
\label{Eqn:NavierStokesEquation}
\partial_t \bm{u} + \bm{u} \cdot \nabla \bm{u}
& = \frac{1}{Re} \triangle \bm{u} - \nabla p& \qquad & \textrm{in} \ \Omega.
\end{alignat}
\end{subequations}
The boundary conditions for the velocity comprise the no-slip
condition at the cylinder, the free-stream condition at the inlet
boundary, the free-slip condition at the lateral boundaries and the
no-stress condition at the downstream outflow boundary. The initial
condition for the velocity at $t=0$ is the unstable fixed point
$\bm{u}_s$ of Navier-Stokes equation \eqref{Eqn:IncompressibleFlow}
perturbed by the real part of the corresponding most unstable
stability eigenmode $\bm{u}_1^{\star}(\bm{x})$
\begin{equation}
\label{Eqn:InitialCondition}
\bm{u} (\bm{x},0) = \bm{u}_s (\bm{x}) + 0.02 \> \bm{u}_1^{\star}(\bm{x}).
\end{equation}
The field $\bm{u}_1^{\star}(\bm{x})$ is normalized to have unit
{$L_2(\Omega)$} norm.

\begin{figure}
\label{Fig:Flows}
\begin{center}
(a) \\
\includegraphics[width=70mm]{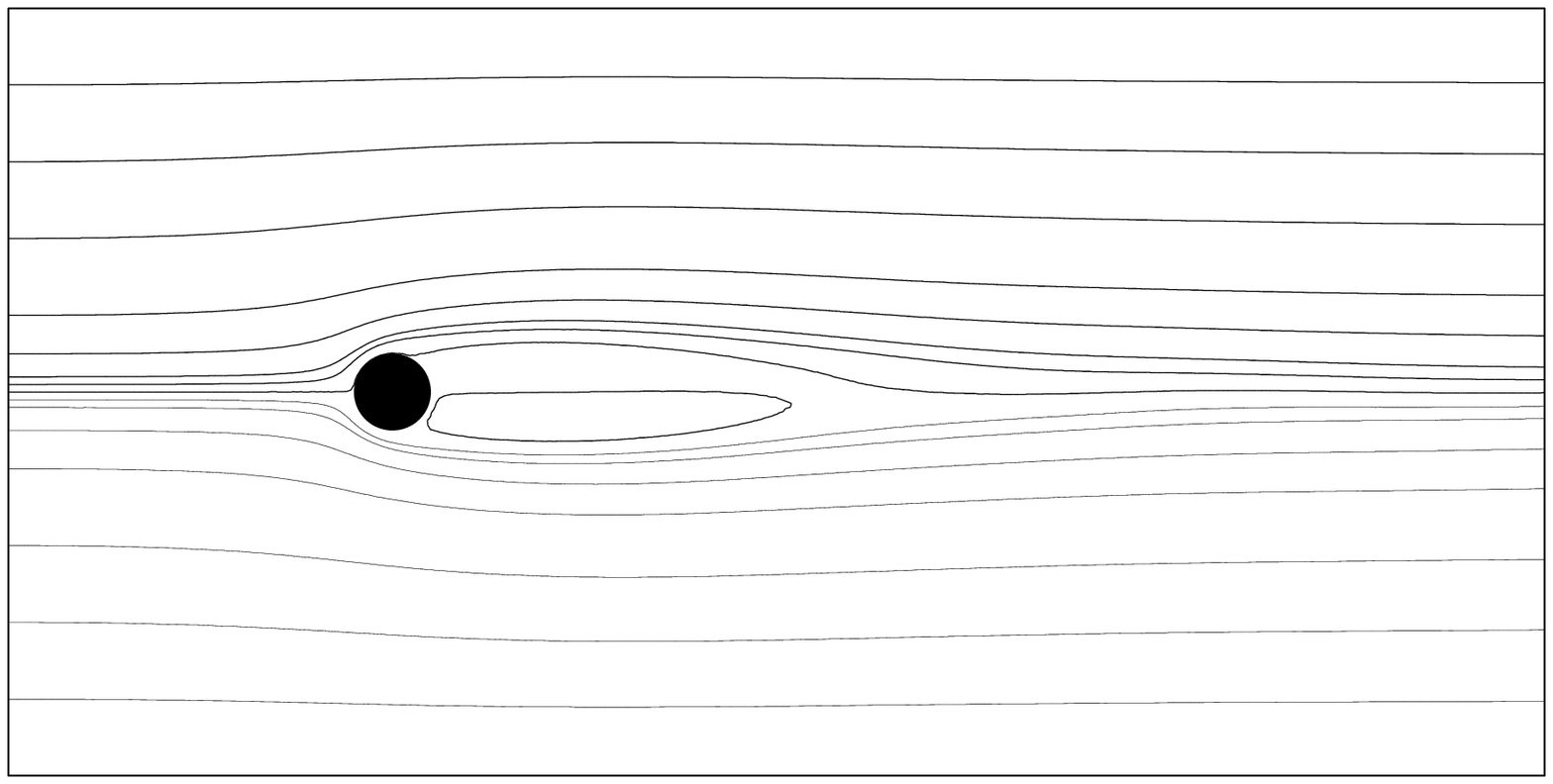} \\
(b) \\
\includegraphics[width=70mm]{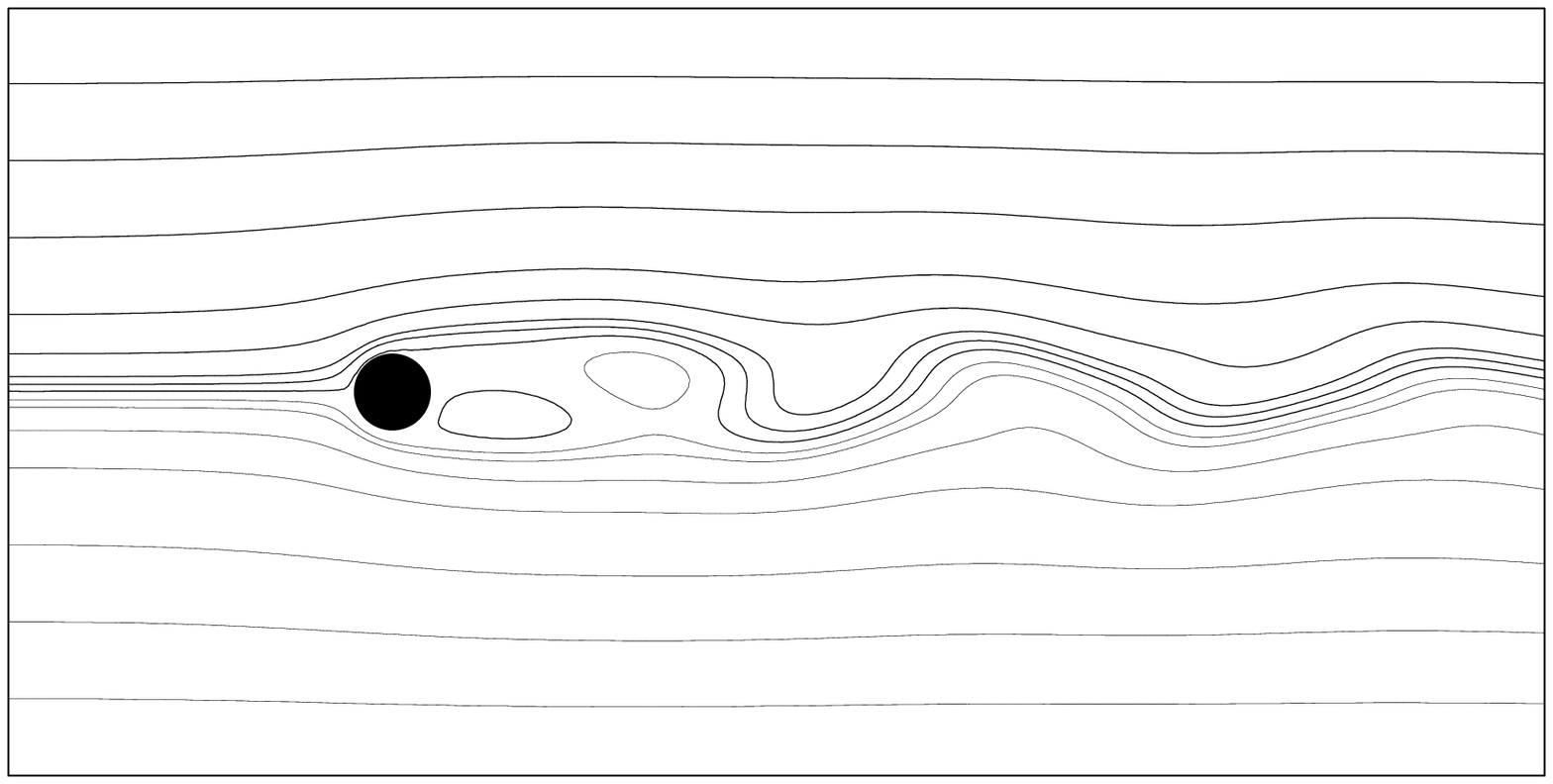} \\
(c) \\
\includegraphics[width=70mm]{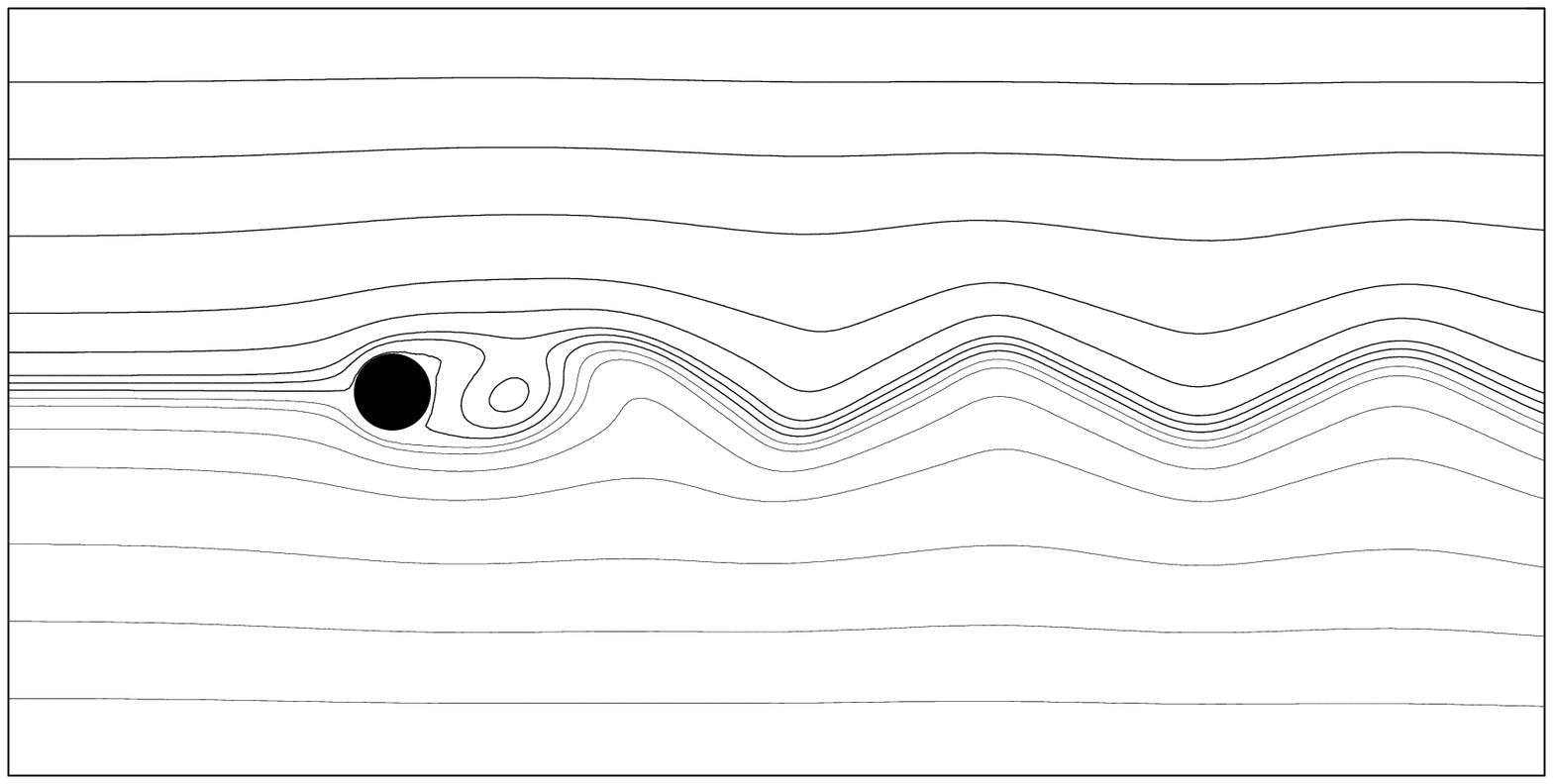} 
\end{center}
\caption[X]{Flow snapshots at (a) $t=0$ (initial condition), (b)
  $t=24.8$ (intermediate transient state), and (c) $t=99.6$ (periodic
  solution at the limit cycle). The domain shown represents the
  computational domain $\Omega$ which also coincides with the region
  where the Galerkin expansion is defined.  The flow patterns are
  visualized using streamlines (which are defined as the level sets of
  the streamfunction $\psi \; : \; \Omega \rightarrow \RR$ related to
  the velocity components via $u = \Dpartial{\psi}{y}$ and $v =
  -\Dpartial{\psi}{x}$).  The cylinder is indicated by the black
  circle.  }
\end{figure}

The initial-boundary-value problem is numerically integrated with a
finite-element method based on an unstructured grid and details are
provided in \cite{Noack2003jfm}.  Figure \ref{Fig:Flows} depicts three
flow snapshots corresponding to the initial condition, an intermediate
transient state and the flow {approaching} the limit cycle.
%----- Figure ----------------------------------------------------------

%=======================================================================
\subsection{Galerkin Expansion}
\label{ToC:GalerkinExpansion}

Below we review approximate modelling approaches typically employed to
obtain low-dimensional descriptions of the dynamics described by
\eqref{Eqn:IncompressibleFlow}, cf.~\cite{hlb96}. The quantities
corresponding to the limit cycle will be denoted with the superscript
'$\circ$' and in the present problem in which the limit cycle is
stable we will therefore have $r^\circ = r_{max}$.  The transient wake
is usually characterized by the superposition of a time-varying
symmetric base flow $\bm{u}^B$ with the length of the recirculating
region decreasing with time as the vortex shedding develops and an
antisymmetric oscillating field $\bm{u}^{\prime}$ representing the
vortex shedding.  The change of the base flow is accurately captured
by a single ``shift mode'' $\bm{u}_{\Delta}$ \cite{Tadmor2010pf},
while the oscillation field can be approximated by the first two
Proper Orthogonal Decomposition (POD) modes $\bm{u}_i$, $i=1,2$,
computed at the limit cycle \cite{Deane1991pfa,Noack2003jfm}.
{Hereafter, the subscript '$\Delta$' will denote various
  quantities related to the shift mode $\bm{u}_{\Delta}$.} Some
information about the construction of orthogonal bases using the POD
approach is presented in Appendix \ref{ToC:POD}. The resulting
truncated Galerkin expansion thus takes the form
\begin{subequations}
\label{Eqn:MeanFieldExpansion}
\begin{align}
\label{Eqn:MeanFieldExpansion_1}
\bm{u}         (\bm{x},t) &\approx \bm{u}^B (\bm{x},t) + \bm{u}^{\prime} (\bm{x},t), \\
\label{Eqn:MeanFieldExpansion_2}
\bm{u}^B       (\bm{x},t) &= \bm{u}_s (\bm{x}) + a_{\Delta}(t) \> \bm{u}_{\Delta} (\bm{x}) ,\\
\label{Eqn:MeanFieldExpansion_3}
\bm{u}^{\prime}(\bm{x},t) &= a_1(t) \> \bm{u}_1 (\bm{x}) + a_2(t) \> \bm{u}_2 (\bm{x}).
\end{align}
\end{subequations}
One typically ignores deformations of the oscillatory modes during the
transient.  These deformations do not significantly alter the model
predictions and do not have any effect on the proposed model
identification approach. On the other hand, inclusion of this effect
in our model would significantly complicate the study and is outside
the scope of this investigation. Identifying $a_3 = a_{\Delta}$ and
assuming phase-invariant behavior, we note that the evolution of the
Galerkin expansion coefficients is governed by a system in the form
\eqref{eq:DescriptorSystemA}, provided that $g_1(0)>0$, corresponding
to an unstable fixed point at the origin, and $g_1(r^{\circ}) = 0$ and
$dg_1/dr|_{r=r^{\circ}} < 0$, corresponding to a locally attracting
limit cycle at $r = r^{\circ} > 0$, cf.~Assumption \ref{ass1}(a).
Moreover, from the assumptions made in Section \ref{ToC:Introduction},
it follows that $\forall_{r \ge 0} \ g_2(r) > 0$. Introducing a number
of further simplifications one arrives at two well-known reduced
models, the Mean-Field Model and the Landau Model.  Since they provide
a point of reference for our test problems, we describe them briefly
below.

%=======================================================================

\subsection{Mean-Field Model}
\label{ToC:MeanFieldSystem}
The mean-field model of oscillatory flow instabilities
\cite{Stuart1971arfm} explicates the mechanism of amplitude saturation
in the mean-field deformation.  This deformation is quantified by the
amplitude $a_{\Delta}$ {of the shift mode,
  cf.~\eqref{Eqn:MeanFieldExpansion_2},} which is slaved to the
fluctuation level $r$.  For the soft (supercritical) Hopf bifurcation,
the resulting evolution equations read
%----- Equation --------------------------------------------------------
\begin{subequations}
\label{Eqn:MeanFieldSystem}
\begin{eqnarray}
\label{Eqn:MeanFieldSystem_1}
\dot r{(t)}      &=& \left[ 
     \sigma_1 - \beta_{\Delta} \> a_{\Delta}({t})
                \right] \> r{(t)},
\\
\label{Eqn:MeanFieldSystem_2}
\dot \theta{(t)} &=& 
     \omega_1 + \gamma_{\Delta} a_{\Delta}({t}),
\\
\label{Eqn:MeanFieldSystem_3}
a_{\Delta}({t}) &=& \alpha_{\Delta} \> r^2{(t)}.
\end{eqnarray}
\end{subequations}
% BRN: I removed the superscript 'M' because it is not explained
% BRN: and the explanation would be 'too heavy' for the normal reader.
%-----------------------------------------------------------------------
The ordinary differential equations \eqref{Eqn:MeanFieldSystem_1}
{and} \eqref{Eqn:MeanFieldSystem_2} correspond to the
Navier-Stokes equation linearized around the time-varying base flow
with $a_{\Delta}$ as the order parameter.  The third algebraic
equation \eqref{Eqn:MeanFieldSystem_3} represents the
Reynolds equation linking the mean flow to the Reynolds stress
{generated by} the fluctuations described {by the} first
equation.

Evidently, equations \eqref{Eqn:MeanFieldSystem} are a particular case
of \eqref{eq:DescriptorSystemA} with $g_1(r) = \sigma_1 -
\beta_{\Delta} \> {a_{\Delta} \> r^2}$, $g_2(r) = \omega_1 + \gamma_{\Delta} \>
{a_{\Delta} \> r^2}$ and $g_3(r) = \alpha_{\Delta} \> r^2$. The parameters
$\sigma_1$, $\omega_1$, $\beta_{\Delta}$, $\gamma_{\Delta}$,
$\alpha_{\Delta}$ of \eqref{Eqn:MeanFieldSystem} may be derived from
the Galerkin approximation described in Section
\ref{ToC:GalerkinExpansion}, or identified from the solutions of the
Navier-Stokes equation via suitable fitting. The periodic solution of
\eqref{Eqn:MeanFieldSystem} is given by
%----- Equation --------------------------------------------------------
\begin{subequations}
\label{Eqn:MeanFieldSystemSolution}
\begin{eqnarray}
r^{\circ}               &=& \sqrt{\frac{\sigma_1}{\alpha_{\Delta} \> \beta_{\Delta}}},
\\
\omega^{\circ}-\omega_1 &=& \gamma_{\Delta} \sqrt{\frac{\sigma_1}{\beta_{\Delta}}},
\\
a_{\Delta}^{\circ} &=& \sigma_1/\beta_{\Delta}.
\end{eqnarray}
\end{subequations}
{At the limit cycle} the mean flow is predicted to have a
vanishing growth rate $g_1$.  This marginal stability property of mean
flows has been conjectured by Malkus \cite{Malkus1956jfm} and is
corroborated by the global stability analysis of the Navier-Stokes
equation \cite{Barkley2006ep}.  Derivation details and a further
discussion of the properties of the mean-field model can be found in
\cite{Noack2003jfm}.

The Landau equation for the supercritical Hopf bifurcation is a
corollary to the mean-field model and constitutes a prototype
evolution equation for self-amplified, amplitude-limited oscillations.
{It is} obtained by substituting \eqref{Eqn:MeanFieldSystem_3} in
\eqref{Eqn:MeanFieldSystem_1}--\eqref{Eqn:MeanFieldSystem_2}:
%----- Equation --------------------------------------------------------
\begin{subequations}
\label{Eqn:LandauEquation}
\begin{eqnarray}
\dot r      &=& \sigma_1 \> r - \beta \> r^3, 
\label{Eqn:LandauEquation_1}
\\
\dot \theta &=& \omega_1      + \gamma \> r^2.
\label{Eqn:LandauEquation_2}
\end{eqnarray}
\end{subequations}
%-----------------------------------------------------------------------
Here, {we require} $\sigma_1,\beta,\omega_1>0$ to ensure a stable
limit cycle with a positive angular velocity in the $(a_1,a_2)$ plane,
cf.~Assumption \ref{ass1}(a), whereas $\gamma$ may vanish, be positive
or negative.  Evidently, equations \eqref{Eqn:LandauEquation}
correspond to
\eqref{Eqn:MeanFieldSystem_1}--\eqref{Eqn:MeanFieldSystem_2} with
$\beta = \alpha_{\Delta} \beta_{\Delta}$ and $\gamma = \alpha_{\Delta}
\gamma_{\Delta}$. {Therefore, the Landau equation
  \eqref{Eqn:LandauEquation} is also easily recognized as a particular
  case of
  \eqref{eq:DescriptorSystemB_1}--\eqref{eq:DescriptorSystemB_2} and
  \eqref{eq:DescriptorSystemA_1}, while the mean-field equation
  \eqref{Eqn:MeanFieldSystem_3} is an example of
  \eqref{eq:DescriptorSystemB_3} and \eqref{eq:DescriptorSystemA_2}.
  The dependency of the growth rate and frequency on the shift mode
  amplitude $a_{\Delta}$ in
  \eqref{Eqn:MeanFieldSystem_1}--\eqref{Eqn:MeanFieldSystem_2} may in
  principle be recovered from
  \eqref{eq:DescriptorSystemB_1}--\eqref{eq:DescriptorSystemB_2} by
  inverting \eqref{eq:DescriptorSystemB_3} to give $r=r(a_{\Delta})$
  and substituting $r(a_{\Delta})$ into
  \eqref{eq:DescriptorSystemB_1}--\eqref{eq:DescriptorSystemB_2}.
  This inversion assumes a monotonous dependence of $a_{\Delta}$ on
  $r$ which is observed in actual wake data.  For completeness, we
  note that the least-order Galerkin model introduced in
  \cite{Noack2003jfm} includes a fast-dynamics equation for
  $a_{\Delta}$.  This dynamic equation is well represented by an
  inertial manifold $a_{\Delta}=a_{\Delta}(r)$ of the form
  \eqref{eq:DescriptorSystemB_3} or \eqref{eq:DescriptorSystemA_2}
  which ignores short transients.}

The Landau equation \eqref{Eqn:LandauEquation} may also be derived
from the first principles, e.g., using a center-manifold reduction
method, see, e.g.,
\cite{Haken1983book,Guckenheimer1986book,Gorban2005book}.  The
parameters may be identified from solutions of the Navier-Stokes
equation with $\sigma_1$ and $\omega_1$ obtained as the real and
imaginary part of the most unstable eigenvalue associated with the
Hopf bifurcation, whereas the nonlinearity parameters $\beta$ and
$\gamma$ can be inferred from the post-transient amplitude and
frequency, $r^{\circ}$ and $\omega^{\circ}$, respectively. In Section
\ref{ToC:Results} we will demonstrate using our proposed approach how
the structure of Landau model \eqref{Eqn:LandauEquation} could be
modified to better reproduce the actual behavior.

%=======================================================================

\subsection{Measurements}
\label{ToC:Measurements}
The state in our approximate model \eqref{eq:DescriptorSystemB} is
characterized by three time-dependent mode amplitudes: $a_1$, $a_2$
and {$a_3$} and as ``measurements'' we will consider the functions
$\ta_1(t)$, $\ta_2(t)$ and $\ta_{\Delta}(t)$ which are obtained for $t
\in [0,T]$ by solving initial-boundary-value problem
\eqref{Eqn:IncompressibleFlow}-\eqref{Eqn:InitialCondition} followed
by {projection, in terms of the inner product defining the POD
  analysis (cf.~Appendix \ref{ToC:POD}), of the resulting
  time-dependent velocity field $\bm{u}(t,\cdot)$ on the modes
  $\bm{u}_1$, $\bm{u}_2$ and $\bm{u}_{\Delta}$. The same procedure
  applies when the velocity field comes from time- and space-resolved
  measurements. In either case, determination of the shift mode
  $\bm{u}_{\Delta}$ requires access to the unstable equilibrium
  $\bm{u}_s$, cf.~\eqref{Eqn:MeanFieldExpansion_2}, which typically
  needs to be obtained numerically.} {The mean-field model may
  also be constructed directly from experimental measurements.  Let
  $s(t)$ be, for instance, a pressure or hot-wire signal with a
  dominant harmonic component corresponding to the vortex shedding.
  Let $\langle s \rangle(t)$ be the short-time mean value, e.g., a
  one-period average, and $s^\prime(t) := s(t) - \langle s \rangle(t)$
  be the fluctuation.  Next, let $a_1$ and $a_2$ be the local cosine
  and sine component of the fluctuation.  These variables may be
  obtained from a Hilbert transform or, more robustly, from a Morlet
  transform of the data.  Finally, we identify $a_{\Delta} = \langle s
  \rangle + \hbox{const}$, where the tunable constant shall ensure the
  correct fixed-point behavior, i.e., $a_{\Delta} =0$ when $r=0$.
  Then, $a_1$, $a_2$ and $a_{\Delta}$ may be approximated by the
  mean-field model. In this case, it is not required that the fixed
  point be actually reached for the identification of the model.  We
  only need a transient with a range of $r$ values over which the
  functions $g_i$ are identified.} In the next Section we introduce a
computational approach allowing one to optimally identify the
functions $g_1(r)$, $g_2(r)$ and $g_3(r)$ in a suitable class, so that
the predictions of system \eqref{eq:DescriptorSystemB} best match in
the least-squares sense the {data} $\ta_1$, $\ta_2$ and
$\ta_{\Delta}$.  We note that, alternatively, a system of evolution
equations for $a_1$, $a_2$ and $a_3$ may be obtained by substituting
ansatz \eqref{Eqn:MeanFieldExpansion_1} into momentum equation
\eqref{Eqn:NavierStokesEquation}, and some comparisons between such an
approach and the results obtained using the model identification
method developed here will be drawn in Section
\ref{ToC:Discussion_physical}.

%***********************************************************************
\section{Computational Approach}
\label{ToC:Computation}

The task of identifying the functions $g_i$, $i=1,2,3$, such that the
output of system \eqref{eq:DescriptorSystemB} matches certain
``measurements'' is an example of an inverse problem \cite{t05}. What
makes this problem somewhat different from typical inverse problems is
that the functions sought have the form of ``constitutive relations'',
in the sense that they depend on the {\em state variables} (i.e., the
dependent variables in the problem), rather than the independent
variables. {More specifically, in the problem considered here} $g_i$,
$i=1,2,3$, depend on $r=\sqrt{a_1^2+a_2^2}$ as opposed to $t$.
Non-parametric formulations of such inverse problems have received
some attention in the context of systems described by PDEs
\cite{bvp10,bp11a,b12,chl74}, but we are not aware of similar
approaches applied to the state-space description of dynamical
systems.

\subsection{Formulation of Optimization Problem}
\label{ToC:Formulation}

We will look for functions $g_i$, $i=1,2,3$, as elements of the
Sobolev space $H^1(\I)$ of continuous functions with square-integrable
gradients on $\I$ which is the ``identifiability'' region defined in
Section \ref{ToC:Introduction}, i.e., the interval spanned by the
state variable $r(t)$ during the system evolution, see also
\cite{bvp10}. Some additional remarks concerning the regularity of
functions $g_1$ and $g_2$ are presented in Appendix
\ref{ToC:Regularity}. In the reconstruction problems considered in
this study, the boundary behavior of the reconstructed functions will
have to be restricted in different ways. This reflects the fact that
system \eqref{eq:DescriptorSystemB} with the reconstructed functions
$g_i$, $i=1,2,3$, should exhibit the behavior dictated by Assumption
\ref{ass1}(a) in specific regions of the phase space, namely, at the
equilibrium $r=0$ and at the limit cycle $r=r^\circ$.  There is some
flexibility as regards possible choices and the boundary conditions we
adopt will reproduce the behavior described by the mean-field model
introduced in Section \ref{ToC:MeanFieldSystem}.  More specifically,
\begin{itemize}
\item
at the origin $r=0$
\begin{itemize}
\item the Jacobian of the right-hand side (RHS) in equation
  \eqref{eq:DescriptorSystemB_1} should be given by $g_1(0)$ which is
  a priori unknown and is to be determined as a part of the solution
  of the reconstruction problem; on the other hand, the Jacobian of
  the RHS in equation \eqref{eq:DescriptorSystemB_2} should vanish;
  this is achieved when the derivatives of both functions $g_1$ and
  $g_2$ are set to zero at the origin, i.e.,
\begin{equation}
\frac{d}{dr}g_i(r)\big|_{r=0}=0, \qquad i=1,2,
\label{eq:BC1}
\end{equation}
\end{itemize}

\item
at the limit cycle $r = r^{\circ}$
\begin{itemize}
\item the RHS of \eqref{eq:DescriptorSystemB_1} should vanish
  resulting in the vanishing of $g_1$, i.e.,
\begin{equation}
g_1(r^\circ) = 0,
\label{eq:BC2}
\end{equation}

\item in regard to equation \eqref{eq:DescriptorSystemB_2}, we will
  prescribe a given slope $G>0$ of the RHS, i.e.,
\begin{equation}
\frac{d}{dr}g_2(r)\big|_{r=r^\circ}=G.
\label{eq:BC3}
\end{equation}
\end{itemize}
\end{itemize}
We note that, while the boundary behavior described by
\eqref{eq:BC1}--\eqref{eq:BC3} at the equilibrium and the limit cycle
is the same as in the mean-field model (cf.~Section
\ref{ToC:MeanFieldSystem}), the behavior of the constitutive relations
$g_1$ and $g_2$ for intermediate values of the state magnitude $0 < r
< r^\circ$ can be arbitrary and will be determined using our optimal
reconstruction procedure. No restrictions are placed on the boundary
behavior of function $g_3$, cf.~\eqref{eq:DescriptorSystemB_3}.

A convenient way to solve inverse problems is to formulate them as
suitable optimization problems \cite{v02}. For each of the functions
$g_i$, $i=1,2,3$, we thus define the corresponding cost functional
$\J_i(g_i) \; : \; H^1(\I) \rightarrow \RR$ as
\begin{subequations}
\label{eq:J}
\begin{align}
\J_1(g_1) & := \frac{1}{2} \int_0^T \left[ r(t) - \tr(t) \right]^2 \, dt,
\label{eq:J1} \\
\J_2(g_2) & := \frac{1}{2} \int_0^T \left[ e^{\imath \theta(t)} - e^{\imath \tilde{\theta}(t)} \right]^2 \, dt,
\label{eq:J2} \\
\J_3(g_3) & := \frac{1}{2} \int_0^T \left[ g_3(r(t)) - \ta_{\Delta}(t) \right]^2 \, dt,
\label{eq:J3} 
\end{align}
\end{subequations}
where $\tr(t) := \sqrt{ \ta_1(t)^2 + \ta_2(t)^2}$, $\tilde{\theta}(t)
:= \arctan( \ta_2(t)/\ta_1(t))$ and $\ta_{\Delta}(t)$ are the
``measurements'' obtained as described in Section
\ref{ToC:Measurements}. {The length $T$ of the assimilation
  window will be chosen sufficiently long to allow the transient
  trajectory to settle on the limit cycle.}  In \eqref{eq:J1} and
\eqref{eq:J2} the functions $r(t)$ and $\theta(t)$ are related to
$g_1$ and $g_2$ via system \eqref{eq:DescriptorSystemB}.  The optimal
reconstructions $\hat{g}_1$, $\hat{g}_2$, and $\hat{g}_3$ are defined
as solutions of the following problems
\begin{subequations}
\label{eq:P}
\begin{alignat}{2}
& (P1) \qquad & \hat{g}_1 & 
:= \argmin_{g_1 \in H^1(\I), \ \frac{d}{dr}g_1(r)|_{r=0}=0, \ g_1(r^\circ) = 0} \ 
\J_1(g_1), \label{eq:P1} \\
& (P2) & \hat{g}_2 & 
:= \argmin_{g_2 \in H^1(\I), \ \frac{d}{dr}g_2(r)|_{r=0}=0,  \ \frac{d}{dr}g_2(r)|_{r=r^\circ}=G} \ \J_2(g_2), \label{eq:P2} \\
& (P3) & \hat{g}_3 & := \argmin_{g_3 \in H^1(\I)} \ \J_3(g_3). \label{eq:P3} 
\end{alignat}
\end{subequations}
Below we describe in detail a gradient-based approach to solution of
Problem $(P1)$. Problem $(P2)$ has a similar structure and the
solution method is essentially the same with some modifications
described hereafter. In Problems $P1$ and $P2$ it is assumed that the
functions, respectively, $g_2$ and $g_1$ are fixed. Problem $(P3)$ has
a different structure warranting a separate solution approach which
will be described further below.

\subsection{Gradient-Based Approach to Solution of Problems $(P1)$ and $(P2)$}
\label{ToC:Grad}

\begin{figure}
\begin{center}
\mbox{
\subfigure[]{
\psfrag{C}[][][0.8]{$\C$}
\includegraphics[width=0.5\textwidth]{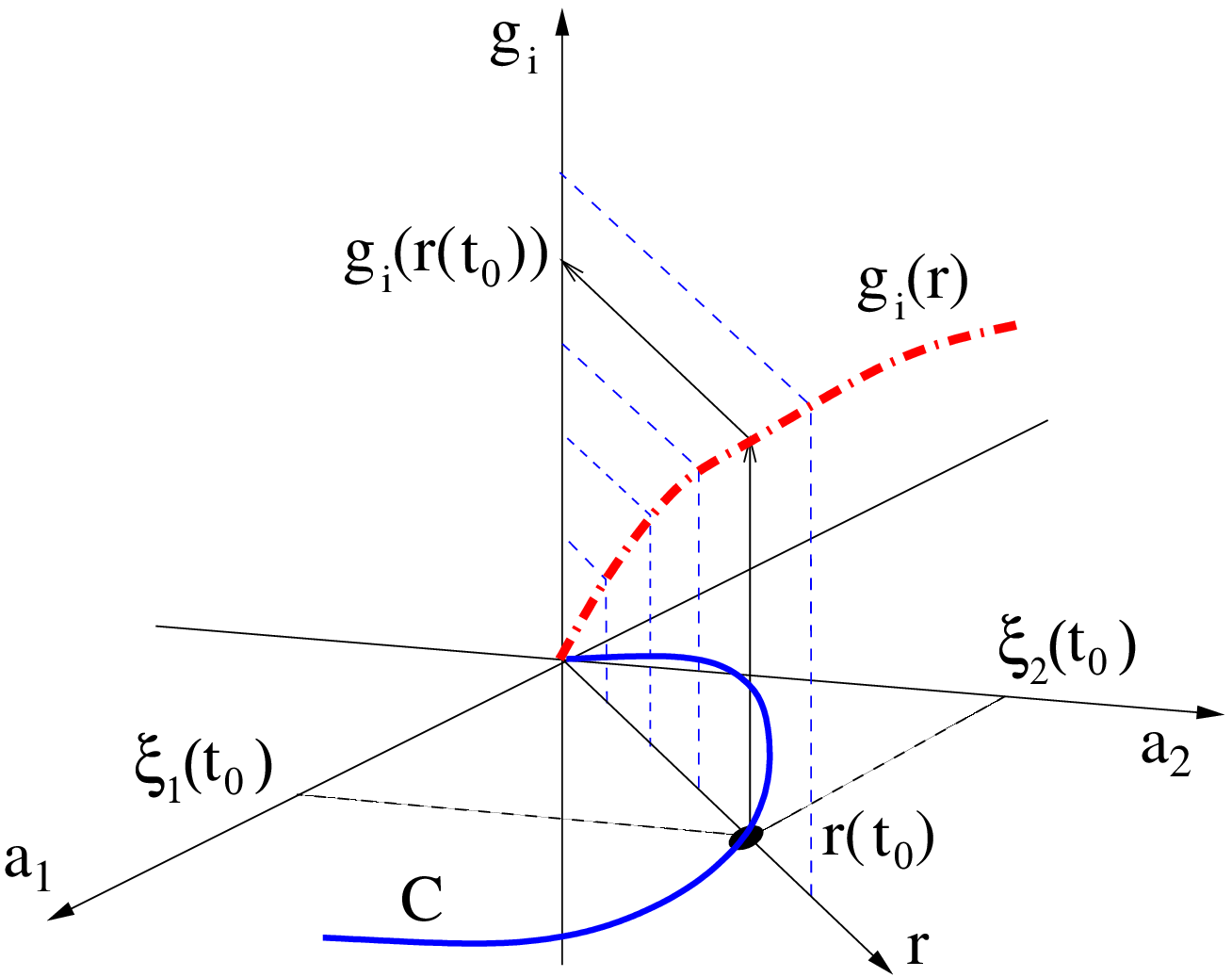}}\qquad
\psfrag{xi}[][][0.8]{$\bxi(t), \, \bxi^*(t)$}
\psfrag{DJ}[][][0.8]{$\nabla \J_i(r)$}
\psfrag{r}[][][0.8]{$r$}
\psfrag{C}[][][0.8]{$\C$}
\subfigure[]{\includegraphics[width=0.4\textwidth]{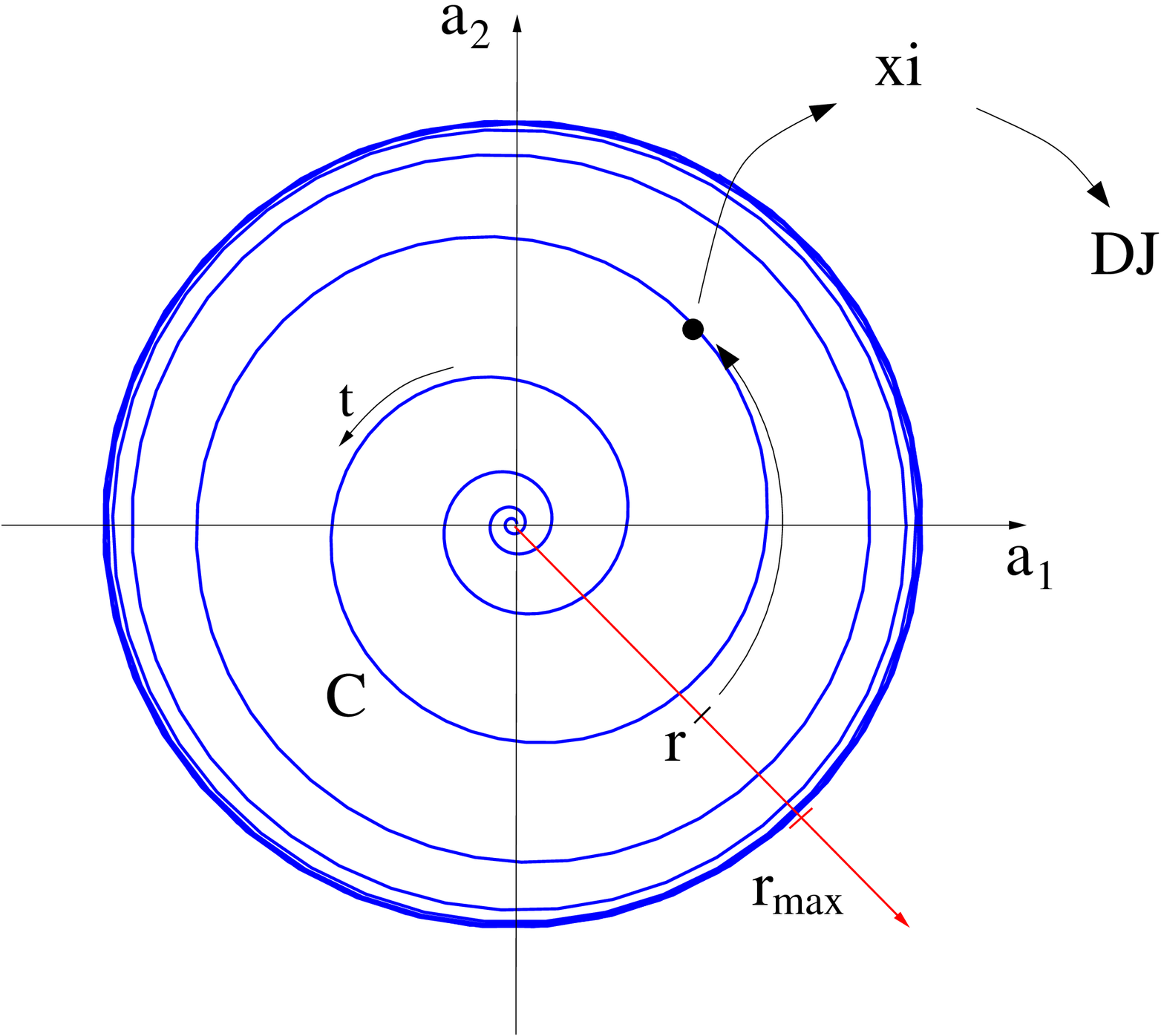}}}
\caption{(a) Schematic indicating the dependence of the ``constitutive
  relation'' on the state magnitude $r$ with the plane $(a_1,a_2)$
  representing the phase space. {Thick red dashed-dotted line:
    function $g_i$, $i=1,2,3$. Thick blue solid line: sample
    trajectory $[\xi_1(t), \ \xi_2(t)]^T$, $t\in [0,T]$, of system
    \eqref{eq:DescriptorSystemA_1}. Black circle: state at time
    $t_0$.} (b) Schematic illustrating the relation between the
  integration variables $dt$ and $dr$, cf.~\eqref{eq:dtdr}. {Blue
    solid line: system trajectory $\C$.} The state $\bxi(t)$ (marked
  with a black point on the trajectory) together with the
  corresponding adjoint state $\bxi^*(t)$ carry information necessary
  to evaluate cost functional gradient $\nabla \J_i(r)$, where $r =
  \|\bxi(t)\|$, cf.~\eqref{eq:gradJL2}.}
\label{fig:vars}
\end{center}
\end{figure}

The minimizer $\hat{g}_1$ is characterized by the first-order
optimality condition \cite{l69} requiring the vanishing of the
G\^ateaux differential \\ $\J_1'(g_1;g'_1) := \lim_{\epsilon\rightarrow
  0} \epsilon^{-1} \left[ \J_1(g_1 + \epsilon g'_1) -
  \J_1(g_1)\right]$, i.e.,
\begin{equation}
\forall_{g'_1 \in H^1(\I), {\ \frac{d}{dr}g'_1(r)|_{r=0}=0, \ g'_1(r^\circ) = 0}} \qquad \J_1'(\hat{g}_1;g'_1) = 0,
\label{eq:opt}
\end{equation}
where $g'_1$ is an arbitrary perturbation direction. The (local)
minimizer can be computed with the following iterative procedure
\begin{equation}
\left\{
\begin{alignedat}{2}
&g_1^{(n+1)} && = g_1^{(n)} - \tau^{(n)} \nabla\J_1(g_1^{(n)}),
 \qquad n=1,\dots, \\
&g_1^{(1)}   && = g_1^0,
\end{alignedat}
\right.
\label{eq:desc}
\end{equation}
where $g_1^0$ represents the initial guess, $n$ denotes the iteration
count and $\nabla\J_1 \; : \; \I \rightarrow \RR$ is the gradient of
cost functional $\J_1$. The length $\tau^{(n)}$ of the step is determined by
solving the following line minimization problem
\begin{equation}
\tau^{(n)} = \argmin_{\tau > 0} \J_1\left(g_1^{(n)} - \tau \nabla\J_1(g_1^{(n)})\right)
\label{eq:taumax}
\end{equation}
which can be done efficiently using standard techniques such as
Brent's method \cite{pftv86}.  For the sake of clarity, formulation
\eqref{eq:desc} represents the steepest-descent method, however, in
practice one typically uses more advanced minimization techniques,
such as the conjugate gradient method, or one of the quasi-Newton
techniques \cite{nw00}.  Evidently, the key element of minimization
algorithm \eqref{eq:desc} is the computation of the cost functional
gradient $\nabla\J_1$. It ought to be emphasized that, while the
governing system \eqref{eq:DescriptorSystemA} is finite-dimensional,
the gradient $\nabla\J_1$ is a function of the state magnitude $r$ and
as such represents a continuous (infinite-dimensional) sensitivity of
cost functional $\J_1(g_1)$ to the perturbations $g'_1 = g'_1(r)$. The
fact that the control variable $g_1$, and hence also the gradient
$\nabla \J_1$, are functions of the {\em state} variable, rather than
the independent variable (Figure \ref{fig:vars}a), will result in cost
functional gradients with structure rather different than encountered
in typical optimization problems for differential equations (see
\cite{bp11a} for some related questions arising in PDE optimization
problems).

In order to identify an expression for the gradient $\nabla \J_1$, we
proceed by computing the G\^ateaux differential of the cost functional
\begin{equation}
\J_1(g_1; g'_1) = \int_0^T \left[r - \tr\right] r'(g_1;g'_1)\, dt =
\int_0^T \frac{r - \tr}{r} \bxi^T\bxi' \, dt =
\Big\langle \nabla\J_1(g_1), g'_1 \Big\rangle_{\X(\I)},
\label{eq:dJ}
\end{equation}
where we used the identity $r' = \bxi^T\bxi' / r$ and the last
equality is a consequence of the Riesz representation theorem
\cite{b77} with $\langle\cdot,\cdot\rangle_{\X(\I)}$ denoting an inner
product in the Hilbert space $\X(\I)$ (to be specified later) of
functions defined on $\I$.  The perturbation variable $\bxi'$ is a
solution of the following perturbation problem (see Appendix
\ref{ToC:Perturbation} for a derivation)
\begin{subequations}
\label{eq:dxi}
\begin{align}
\dot{\bxi}'(t) &= \left[ g_1(r(t)) \bI + \bI\, \bxi(t)\, (\bnabla g_1(r(t)))^T \right. \nonumber \\
& \hspace*{1.0cm} \left. + g_2(r(t)) \bJ + \bJ\, \bxi(t)\, (\bnabla g_2(r(t)))^T
\right] \bxi'(t) + \bI \, \bxi(t)\,  g'_1 \nonumber \\
& =: \bA(\bxi(t)) \, \bxi' + \bI \, \bxi(t)\,  g'_1, \label{eq:dxi1} \\
\bxi'(0) &= 0, \label{eq:dxi2}
\end{align}
\end{subequations}
where $\bnabla g_i = \left[\Dpartial{g_i}{a_1}, \
  \Dpartial{g_i}{a_2} \right]^T$, $i=1,2$. We note that G\^ateaux
differential \eqref{eq:dJ} is not yet in the form consistent with the
Riesz representation, since the perturbation $g'_1$ does not appear in
it as a factor, but is hidden on the RHS in perturbation
equation \eqref{eq:dxi1}. A standard technique to convert G\^ateaux
differential \eqref{eq:dJ} to the Riesz form is based on the {\em
  adjoint} variable $\bxi^* \; : \; [0,T]\rightarrow \RR^2$. Taking
the inner product (in $\RR^2$) of $\bxi^*(t)$ with equation
\eqref{eq:dxi1}, integrating over $[0,T]$ and then integrating by
parts we obtain
\begin{equation}
\begin{split}
0 = & \int_0^T \, \left( \bxi^*\right)^T 
\Big\{\dot{\bxi}' - \left[ g_1(r) \bI + \bI\, \bxi\, (\bnabla g_1(r))^T + 
g_2(r) \bJ + \bJ\, \bxi\, (\bnabla g_2(r))^T
\right] \bxi' - \bxi\,  g'_1\Big\}\, dt \\ 
= & \int_0^T \, \left( \bxi'\right)^T 
\Big\{-\dot{\bxi}^* - \left[ g_1(r) \bI + \bI\, \bxi\, (\bnabla g_1(r))^T + 
g_2(r) \bJ + \bJ\, \bxi\, (\bnabla g_2(r))^T \right]^T \bxi^* \Big\} \, dt \\
& + \Big[( \bxi^*)^T \bxi'\Big]_{t=0}^{t=T} 
- \int_0^T \, \left( \bxi^*\right)^T \, \bI \, \bxi(t)\,  g'_1 \, dt.
\end{split}
\label{eq:id1}
\end{equation}
Defining the {\em adjoint} system as
\begin{subequations}
\label{eq:axi}
\begin{align}
- \dot{\bxi}^*(t) &= \left[ \bA(\bxi(t)) \right]^T \bxi^*(t) 
+ \frac{r - \tr}{r} \bxi, \label{eq:axi_1} \\
\bxi^*(T) &= 0, \label{eq:axi_2}
\end{align}
\end{subequations}
we reduce relation \eqref{eq:id1} to
\begin{equation}
\J'_1(g_1;g'_1) =  \int_0^T \, \left( \bxi^*\right)^T \, \bI \, \bxi(t)\,  g'_1 \, dt.
\label{eq:dJb}
\end{equation}
We note that, although $g'_1$ already appears as a factor in
expression \eqref{eq:dJb}, this expression is still not in the Riesz
form, since the integration is with respect to the time $dt$, whereas
in the inner product $\langle\cdot,\cdot\rangle_{\X(\I)}$ defining the
Riesz representer integration is with respect to the measure $dr$
defined on the interval $\I$ (connection between the different
integration variables is illustrated schematically in Figure
\ref{fig:vars}b). The two variables are related via the following
transformation
\begin{equation}
r = | \bxi |  = \sqrt{\xi_1^2 + \xi_2^2} \quad \Longrightarrow \quad 
dr = \frac{\xi_1 d\xi_1 + \xi_2 d\xi_2}{r} = \frac{\xi_1 f_1 + \xi_2 f_2}{r} dt,
\label{eq:dtdr}
\end{equation}
where we used the identities $d\xi_1 = f_1 \, dt$ and $d\xi_2 = f_2 \,
dt$, cf.~\eqref{eq:DescriptorSystemA}. Denoting the
trajectory in the state space $\C := \left\{ \cup_{t \in [0,T]} \
  \bxi(t) \in \RR^2 \right\}$, and combining \eqref{eq:dJb} with
\eqref{eq:dtdr} we obtain the expression 
\begin{equation}
\J'_1(g_1;g'_1) =  \int_{\C} \, \frac{\left( \bxi^*\right)^T \, \bI \, \bxi}{\xi_1 f_1 + \xi_2 f_2}\,  g'_1(r) \, dr 
= \int_0^{r_{max}} \, \frac{\left( \bxi^*\right)^T \, \bI \, \bxi}
{\xi_1 f_1 + \xi_2 f_2}\,  g'_1(r) \, dr
\label{eq:dJc}
\end{equation}
which is already in the required Riesz form. In \eqref{eq:dJc} the
line integral over the contour $\C$ and the definite integral over the
interval $\I$ are equal, because for dynamical system
\eqref{eq:DescriptorSystemA_1} points on the contour $\C$ with the
magnitude $r \in (0,r^\circ)$ are unique, so that the map $r
\rightarrow \left( \bxi(t) \, | \, |\bxi(t)| = r \right) \in \C$ is
one-to-one (in more general situations when this is not the case, or
when the reconstructed function depends on more than one state
variable, e.g., both $\xi_1$ and $\xi_2$ here, the change of variables
needed to obtain the Riesz form will be more complicated and one has
to employ more general techniques such as those developed in
\cite{bvp10,bp11a}).

While this is not the gradient we will use in actual computations, we
will first obtain an expression for the $L_2$ gradient which in the
next subsection will be used as the basis for the calculation of
gradients defined in the Sobolev space $H^1(\I)$.  Thus, setting $\X =
L_2(\I)$ in \eqref{eq:dJ}, we obtain from \eqref{eq:dJc}
\begin{equation}
\nabla^{L_2} \J_1(r) = \frac{\left( \bxi^*\right)^T \, \bI \, \bxi(t)}{\xi_1 f_1 + \xi_2 f_2}, \quad \forall_{r \in \I}.
\label{eq:gradJL2}
\end{equation}
Expression \eqref{eq:gradJL2} is validated computationally in Section
\ref{ToC:Validation}.

As regards Problem $P2$, the optimality condition takes the form,
cf.~\eqref{eq:opt} and \eqref{eq:dJ},
\begin{equation}
\begin{aligned}
\forall_{g'_2 \in H^1(\I), {\ \frac{d}{dr}g'_2(r)|_{r=0}=0,  \ \frac{d}{dr}g'_2(r)|_{r=r^\circ}=G}} \quad 
\J'_2(\hat{g}_2;g'_2) 
&=  \int_0^T \, \sin(\theta - \tilde{\theta}) \theta' \, dt \\
&= \int_0^T \, \sin(\theta - \tilde{\theta}) \frac{\bxi^T \, \bJ \, \bxi'}{r^2} \, dt = 0,
\end{aligned}
\label{eq:dJ2}
\end{equation}
where we used the identity $\theta' = r^{-2}\, \bxi^T \, \bJ \,
\bxi'$. Following the same steps as described above, we obtain an
expression for the cost functional gradient in the form
\eqref{eq:gradJL2}, however, the adjoint system satisfied by $\bxi^*$
has now a different source term on the RHS
\begin{subequations}
\label{eq:axi2}
\begin{align}
- \dot{\bxi}^*(t) &= \left[ \bA(\bxi(t)) \right]^T \bxi^*(t) 
+ \frac{\sin(\theta - \tilde{\theta})}{r^2} \bJ\, \bxi(t), \label{eq:axi2_1} \\
\bxi^*(T) &= 0. \label{eq:axi2_2}
\end{align}
\end{subequations}

In regard to Problem $P3$, using change of variables \eqref{eq:dtdr},
we can rewrite \eqref{eq:J3} as
\begin{equation}
\J_3(g_3) =  \frac{1}{2} \int_{\C} \, \frac{r}{\xi_1 f_1 + \xi_2 f_2} \left[ a_3(r) - g_3(r)\right]^2 \, dr.
\label{eq:J3b}
\end{equation}
Thus, in optimization problem \eqref{eq:P3} we look for a function
$g_3 \in H^1(\I)$ which is as close as possible (in a weighted $L_2$
topology) to a given function $a_3 \in L_2(\I)$, This problem, in
fact, does not have a solution because of the density of the function
space $H^1(\I)$ in $L_2(\I)$, cf.~\cite{af05}. However, it is possible
(and satisfactory from the application point of view) to ``solve''
problem \eqref{eq:P3} approximately by finding a $\hat{g}_3 \in
H^1(\I)$ such that $\J_3(g_3)$ is sufficiently small. Such an approach
is described in Section \ref{ToC:Sobolev}.

\subsection{Sobolev Gradients}
\label{ToC:Sobolev}

In this Section we describe how Sobolev gradients $\nabla^{H^1} \J_i
\in H^1(\I)$, $i=1,2$, used in minimization algorithm \eqref{eq:desc}
for Problems $P1$ and $P2$ can be obtained from \eqref{eq:dJc}. In
addition to enforcing smoothness of the reconstructed functions, this
formulation allows us to impose the desired behavior at the
endpoints of the interval $\I$, cf.~\eqref{eq:BC1}-\eqref{eq:BC3},
via suitable boundary conditions. We begin by defining the $H^1$ inner
product on $\I$ as
\begin{equation}
\forall_{z_1,z_2 \in H^1(\I)} \qquad
\big\langle z_1, z_2 \big\rangle_{{H^1(\I)}} = 
\int_0^{r_{max}} z_1 z_2 + \ell^2 \Dpartial{z_1}{r} \Dpartial{z_2}{r}\, dr,
\label{eq:ipH1}
\end{equation}
where $\ell \in \RR$ is a parameter with the meaning of a ``length
scale''. It is well known \cite{pbh04} that extraction of cost
functional gradients in the space $H^1$ with the inner product defined
as in \eqref{eq:ipH1} can be regarded as low-pass filtering the $L_2$
gradients with the cut-off wavenumber given by $\ell^{-1}$.  As
regards the behavior of the gradients $\nabla^{H^1} \J$ at the
endpoints of the interval $\I$, we can require the vanishing of either
the gradient itself or its derivative $\frac{d}{dr}(\nabla^{H^1} \J)$,
and the boundary conditions we prescribe correspond to relations
\eqref{eq:BC1}--\eqref{eq:BC3} introduced as a part of the formulation
of optimization problems \eqref{eq:P1}--\eqref{eq:P2}, cf.~Assumption
\ref{ass1}(a). As regards the boundary data at $r=r^\circ$ (i.e., at
the limit cycle), in Problem $P2$ we use $\frac{d}{dr} \nabla^{H^1}
\J_2(r)|_{r=r^\circ} = 0$ which ensures that the property
$\frac{d}{dr} g_2^0(r)|_{r=r^\circ} = G$ of the initial guess $g_2^0$
remains unchanged during iterations \eqref{eq:desc}.

Identifying expression \eqref{eq:dJc} with inner product
\eqref{eq:ipH1}, cf.~\eqref{eq:dJ}, integrating by parts and using the
boundary conditions mentioned above we obtain the following elliptic
boundary-value problem {on $\I$} defining the Sobolev gradient
$\nabla^{H^1} \J$
\begin{subequations}
\label{eq:gradJH1}
\begin{alignat}{2}
\left( 1 - \ell^2 \frac{d^2}{dr^2} \right) \nabla^{H^1} \J &= \nabla^{L_2} \J & \qquad & \textrm{in} \ (0,r^\circ), \label{eq:gradJH1_1} \\
\frac{d}{dr} \nabla^{H^1} \J &= 0 & & \textrm{at} \ r = 0, \label{eq:gradJH1_2} \\
\left.
\begin{alignedat}{2}
& (P1): & \qquad\quad &   \nabla^{H^1} \J \\
& (P2): &&   \frac{d}{dr} \nabla^{H^1} \J
\end{alignedat}\right\}  &= 0 & & \textrm{at} \ r = r^{\circ}, \label{eq:gradJH1_3} 
\end{alignat}
\end{subequations}
where the expression for $\nabla^{L_2} \J$ is given in \eqref{eq:gradJL2}.

As concerns Problem $P3$, we propose to reconstruct $\hat{g}_3 \in
H^1(\I)$ directly (i.e., without iterations) by solving the following
problem
\begin{subequations}
\label{eq:gradJ3}
\begin{alignat}{2}
\left( 1 - \ell^2 \frac{d^2}{dr^2} \right) g_3^{\ell} &= a_3(r) & \qquad & \textrm{in} \ (0,r^\circ), \label{eq:gradJ3_1} \\
\frac{d}{dr} g_3^{\ell} &= 0 & & \textrm{at} \ r = 0, \label{eq:gradJ3_2} \\
g_3^{\ell} &= a_3(r^\circ) & & \textrm{at} \ r = r^{\circ}, \label{eq:gradJ3_3} 
\end{alignat}
\end{subequations}
which, except for the boundary condition at $r = r^{\circ}$, has an
identical structure as \eqref{eq:gradJH1}. The superscript in
$g_3^{\ell}$ indicates dependence of the solution on the parameter
$\ell$.  We note that as $\ell \rightarrow 0$ the left-hand side (LHS)
in \eqref{eq:gradJ3_1} approaches the identity transformation which
means that $\| g_3^{\ell} - a_3 \|_{L_2(\I)} \rightarrow 0$, so that
also $\J_3(g_3^{\ell}) \rightarrow 0$, {as $\ell \rightarrow 0$}.
Since solutions of system \eqref{eq:gradJ3} are not defined for $\ell
= 0$, we will obtain our approximate reconstruction as $\hat{g}_3 :=
g_3^{\ell}$ for some small value of $\ell$. Results concerning model
identification for the system described in Section
\ref{ToC:ModelIdentification} are presented in Section
\ref{ToC:Results}, whereas in Sections \ref{ToC:Validation} and
\ref{ToC:Discussion_computational} we analyze certain computational
aspects of the method.
%***********************************************************************
\section{Results}
\label{ToC:Results}

\begin{figure}
\begin{center}
\subfigure[]{\includegraphics[width=0.6\textwidth]{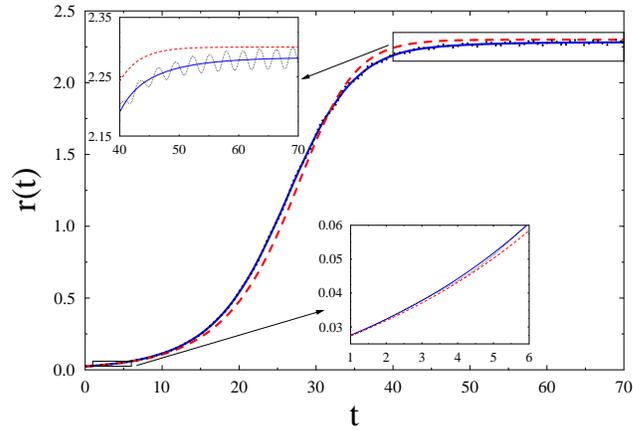}}\vspace*{-0.5cm}
\subfigure[]{\includegraphics[width=0.6\textwidth]{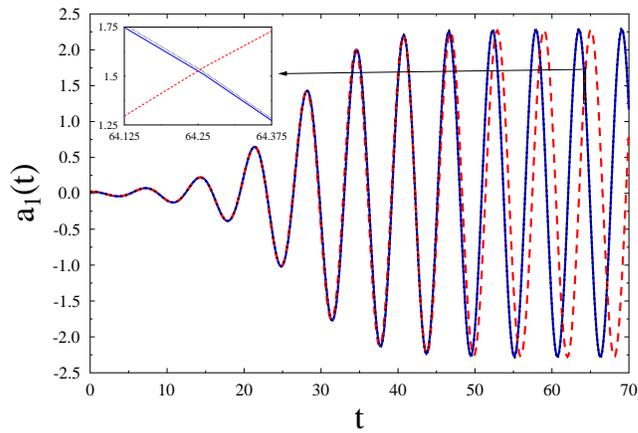}}\vspace*{-0.5cm}
\subfigure[]{\includegraphics[width=0.6\textwidth]{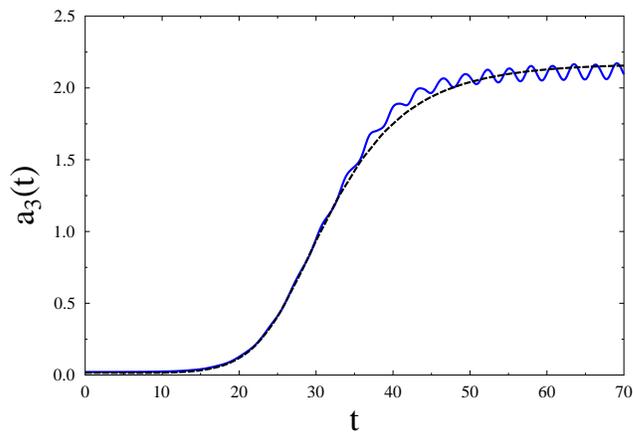}}
\caption{Comparison of (a) the state magnitude $r$, (b) state variable
  $a_1$ and (c) state variable $a_3$ as functions of time.
  {Black dotted lines: measurement data. Blue solid lines:
    solution of system \eqref{eq:DescriptorSystemB} using the optimal
    reconstructions $\hat{g}_1$, $\hat{g}_2$ and $\hat{g}_3$.  Red
    dashed lines: solution of system \eqref{eq:DescriptorSystemB}
    using initial guesses \eqref{eq:IG}.} In Figures (a) and (b)
  insets are included to highlight the differences between the data
  sets.}
\label{fig:ra1a3}
\end{center}
\end{figure}

In this Section we present results concerning the solution of model
identification problems $P1$, $P2$ and $P3$,
cf.~\eqref{eq:P1}--\eqref{eq:P3}, for the system introduced in Section
\ref{ToC:ModelIdentification}. Motivated by practical considerations,
we make the following
\begin{assumption}
\ In the solution of Problems $P1$ and $P2$ we set in system \eqref{eq:DescriptorSystemA_1}
\begin{subequations}
\label{eq:g1g2}
\begin{alignat}{2}
& (P1): & \qquad\qquad  g_2 &= 0, \label{eq:g1g2_1} \\
& (P2): & \qquad\qquad  g_1 &= \hat{g}_1, \label{eq:g1g2_2} 
\end{alignat}
\end{subequations}
\label{ass2}
\end{assumption}
\noindent%
which means that in the reconstruction of $g_2$ we use the best
available estimate of $g_1$ obtained from the solution of Problem
$P1$. The choice of $g_2$ has no effect on the solution of Problem
$P1$ and hence without loss of generality we can adopt
\eqref{eq:g1g2_1}.  As the gradient descent algorithm in Problems $P1$
and $P2$ we use the Polak-Ribiere version of the nonlinear conjugate
gradient method \cite{nw00} in which the ``momentum'' term is reset to
zero every 20 iterations and iterations \eqref{eq:desc} are declared
converged when $(\J_i(g_i^{(n+1)}) - \J_i(g_i^{(n)})) /
\J_i(g_i^{(n)}) \le 10^{-7}$, $i=1,2$. In the system described in
Section \ref{ToC:ModelIdentification} the limit cycle is characterized
by $r^{\circ} = 2.3$, whereas the length of the time window is chosen
as $T = 70$ which is long enough to allow the transient to settle on
the limit cycle (see Figure \ref{fig:ra1a3}a below). We have
successfully solved Problems $P1$, $P2$ and $P3$ using different
combinations of numerical parameters, and the parameters used to
obtain the results presented in this Section are summarized below.
Systems \eqref{eq:DescriptorSystemA} and \eqref{eq:axi} were solved
using MATLAB subroutine {\tt ode45} with an adaptive time-stepping.
Unless stated otherwise, the integrals defined on the interval
$[0,T]$, cf.~\eqref{eq:J}, were discretized using $N_T = 500$
equispaced points. The interval $\I$ was discretized using $N_\I = 75$
equispaced points, and boundary-value problems \eqref{eq:gradJH1} and
\eqref{eq:gradJ3} were approximated using the second-order finite
differences. The length-scale parameter appearing in
\eqref{eq:gradJH1} and \eqref{eq:gradJ3} was $\ell = 1.0$ in Problems
$P1$ and $P2$, and $\ell = 0.1$ in Problem $P3$. The initial condition
$\bxi^0$ for system \eqref{eq:DescriptorSystemA} and the initial
guesses $g_1^0$ and $g_2^0$ for reconstruction algorithm
\eqref{eq:desc} must be chosen so that the magnitudes $|\bxi(t)|$, $t
\in [0,T]$, span the entire interval $\I$, as otherwise the
sensitivities (gradients) cannot be properly defined for {\em all}
values of $r$ (we refer the reader to \cite{bvp10} for a discussion
how this limitation can be overcome in some cases). Since our goal is
now to assess possible improvements to mean-field model
\eqref{Eqn:MeanFieldSystem}, we will use it with the coefficients
determined as discussed in Section \ref{ToC:MeanFieldSystem} as the
initial guess for the reconstructions, so that
\begin{equation}
g_1^0(r) = {0.151  - 0.151} \left(\frac{r}{r^{\circ}}\right)^2, \qquad
g_2^0(r) = 0.886 + 0.15 \left(\frac{r}{r^{\circ}}\right)^2,
\label{eq:IG}
\end{equation}
respectively, for Problems $P1$ and $P2$. It is clear that initial
guesses \eqref{eq:IG} satisfy properties
\eqref{eq:BC1}--\eqref{eq:BC3} with $G=0.224$. As the initial
condition $\bxi^0$ for \eqref{eq:DescriptorSystemA} we used a small
perturbation around the fixed point at the origin.

\begin{figure}
\begin{center}
\subfigure[]{\includegraphics[width=0.6\textwidth]{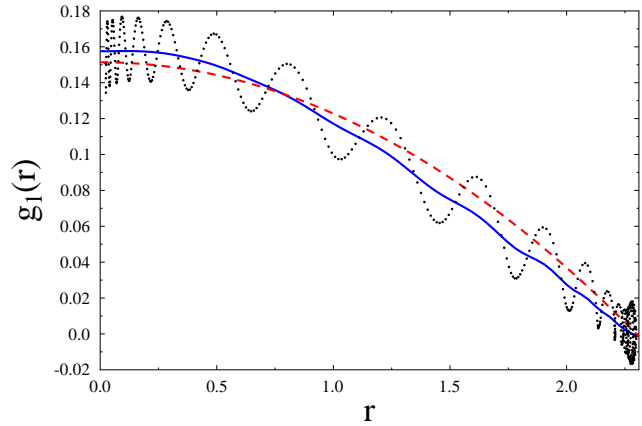}}\vspace*{-0.5cm}
\subfigure[]{\includegraphics[width=0.6\textwidth]{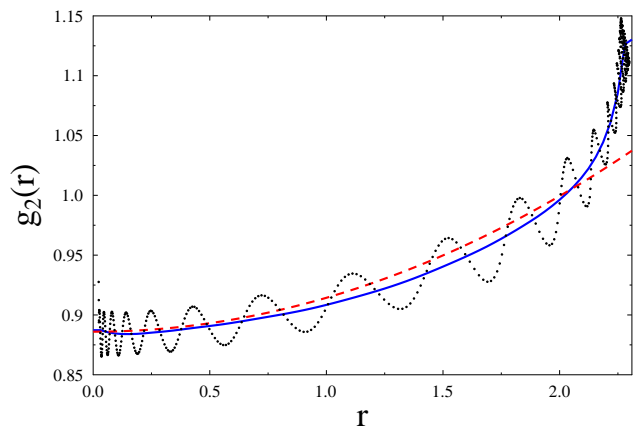}}\vspace*{-0.5cm}
\subfigure[]{\includegraphics[width=0.6\textwidth]{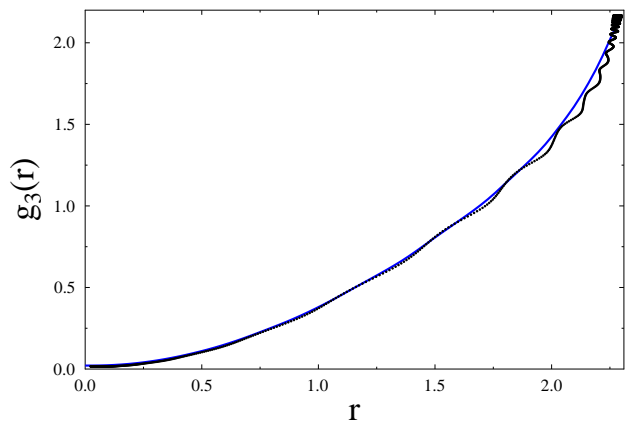}}
\caption{{Blue solid lines:} optimal reconstructions of the
  constitutive relations (a) $\hat{g}_1(r)$, (b) $\hat{g}_2(r)$ and
  (c) $\hat{g}_3(r)$. {Red dashed lines:} the corresponding initial
  guesses (a) $g^{0}_1(r)$ and (b) $g^{0}_2(r)$. {Black symbols:}
  values of (a) $r^{-1} \, (dr/dt)|_{r(t_i)}$, (b)
  $(d\theta/dt)|_{r(t_i)}$ and (c) $a_{\Delta}|_{r(t_i)}$, computed
  based on the measurement data at the time instants $t_i$,
  $i=1,\dots, N_T$.}
\label{fig:g123}
\end{center}
\end{figure}
\begin{figure}
\begin{center}
\includegraphics[width=0.7\textwidth]{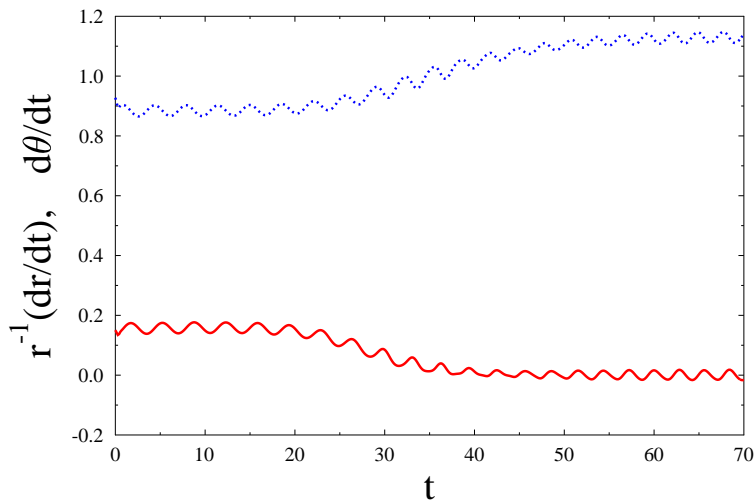}
\caption{Dependence of the LHS in equations
  \eqref{eq:DescriptorSystemB_1}--\eqref{eq:DescriptorSystemB_2}
  {evaluated based on the measurement data on time $t$.  Red
    solid line: $r^{-1} \, (dr/dt)$. Blue dotted line:
    $(d\theta/dt)$.}}
\label{fig:g12_vs_t}
\end{center}
\end{figure}
Our main results are presented in Figure \ref{fig:g123} where we show
the optimal reconstructions $\hat{g}_i(r)$, $i=1,2,3$, and compare
them against the left-hand sides of equations
\eqref{eq:DescriptorSystemB_1}--\eqref{eq:DescriptorSystemB_3}, all
shown as functions of the state magnitude $r$.  For completeness, the
LHS of equations
\eqref{eq:DescriptorSystemB_1}--\eqref{eq:DescriptorSystemB_2} are
shown as functions of time $t \in [0,T]$ in Figure \ref{fig:g12_vs_t}
(in the case of $\hat{g}_1$, cf.~Figures \ref{fig:g123}a and
\ref{fig:g12_vs_t}, the LHS of equation \eqref{eq:DescriptorSystemB_1}
is additionally divided into $r$).  In Figures \ref{fig:g123}a,b we
also indicate the RHS of mean-field model \eqref{eq:IG} which were
used as the initial guesses for the reconstructions. We see in Figures
\ref{fig:g123} that, as expected, the reconstructed constitutive
relations $\hat{g}_i(r)$, $i=1,2,3$, smoothly approximate the
left-hand sides of the corresponding equations evaluated using the
measurements. Considered as functions of $r$, these left-hand sides
are multi-valued which is a consequence of the fact that, due to the
oscillations of the measurement data at the limit cycle (see Figure
\ref{fig:ra1a3}a below), the map $t \rightarrow \tilde{r}(t)$ is not
one-to-one. In Figures \ref{fig:g123}a,b we observe systematic
deviations of the optimal reconstructions $\hat{g}_1$ and $\hat{g}_2$
from the corresponding functions in mean-field model \eqref{eq:IG}. In
addition, based on the reconstructions $\hat{g}_1$ and $\hat{g}_2$ we
can obtain estimates of two important quantities, namely, the growth
rate of the instability at the origin given by $\frac{d}{dr} \left[
  \hat{g}_1 \, r \right]_{r=0} = \hat{g}_1(0) = 0.1576$, and the
oscillation frequency at the limit cycle given by
$\hat{g}_2(r^{\circ}) = 1.130$. These numbers should be compared with,
respectively, 0.151 and 1.036 obtained as discussed in Section
\ref{ToC:MeanFieldSystem} and used in mean-field model \eqref{eq:IG}.
Finally, in Figure \ref{fig:ra1a3} we compare the outputs from system
\eqref{eq:DescriptorSystemB} obtained using the mean-field model and
the reconstructions $\hat{g}_i(r)$, $i=1,2,3$, against the
corresponding measured quantities (as regards the time-history of the
state variables, we do not show $a_2(t)$, as it has qualitatively very
similar behavior to $a_1(t)$ already shown in Figure
\ref{fig:ra1a3}b). In Figures \ref{fig:ra1a3}a,b (see, in particular,
the insets) we note that the evolution of $r(t)$ and $a_1(t)$ obtained
using the optimal reconstructions $\hat{g}_1$ and $\hat{g}_2$ is much
closer to the measured quantities than the evolutions computed using
mean-field model \eqref{eq:IG}.  We remark, however, that the
measurements $\tr(t)$ shown in Figure \ref{fig:ra1a3}a reveal some
high-frequency oscillations which are not captured by the trajectory
$r(t)$ obtained using the optimal reconstruction $\hat{g}_1$. These
oscillations reflect a phase dependence in the behavior of the
solutions of the original Navier-Stokes equation
\eqref{Eqn:IncompressibleFlow}, an effect which is by construction
excluded from ansatz \eqref{eq:DescriptorSystemB_1}, cf.~Assumption
\ref{ass1}(a). As a consequence, $\hat{g}_1$ can convey phase-averaged
information only.  We will return to this problem again in Section
\ref{ToC:Discussion_physical}.  As regards the results shown in Figure
\ref{fig:ra1a3}c, we see that, while the optimal reconstruction
$\hat{g}_3(r)$ is quite smooth (Figure \ref{fig:g123}c), the quantity
$\hat{g}_3(r(t))$ exhibits oscillations absent in the original
measurement data $\ta_3(t)$. This effect as well is a consequence of
the lack of phase-dependence in ansatz \eqref{eq:DescriptorSystemB_1}
and \eqref{eq:DescriptorSystemB_3}.

%***********************************************************************
\section{Discussion of Computational and Physical Aspects}
\label{ToC:Discussion}

In this Section we analyze a number of computational and physical
modelling aspects of the proposed approach which can be important in
applications. We begin by examining the accuracy of the cost
functional gradients in Section \ref{ToC:Validation}, followed by a
study of the robustness of iterations \eqref{eq:desc} in Section
\ref{ToC:Discussion_computational} and conclude with some insights
about physical modelling in Section \ref{ToC:Discussion_physical}.

\subsection{Validation of Gradients}
\label{ToC:Validation}

\begin{figure}
\begin{center}
\mbox{
\subfigure[]{\includegraphics[width=0.5\textwidth]{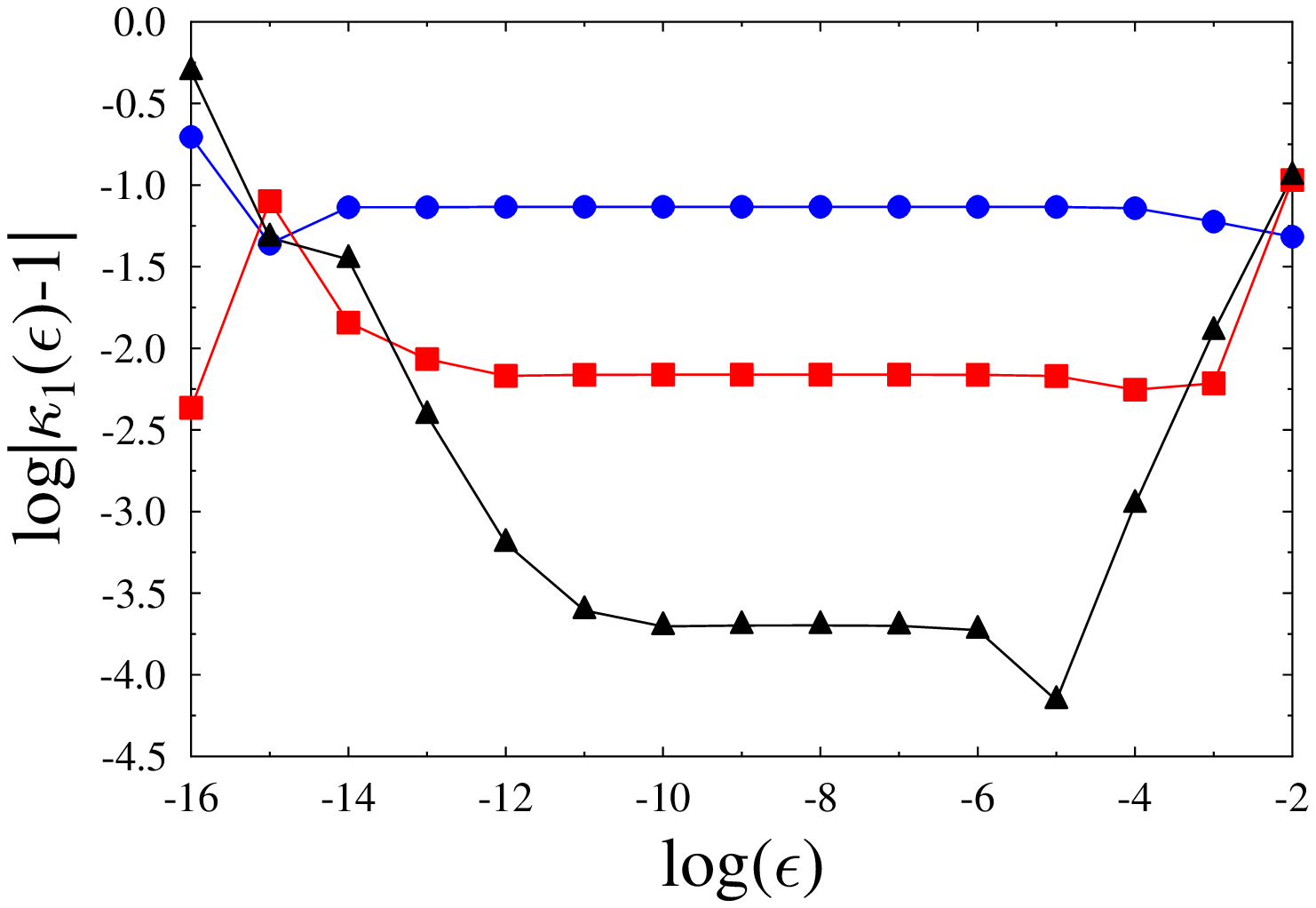}}\quad
\subfigure[]{\includegraphics[width=0.5\textwidth]{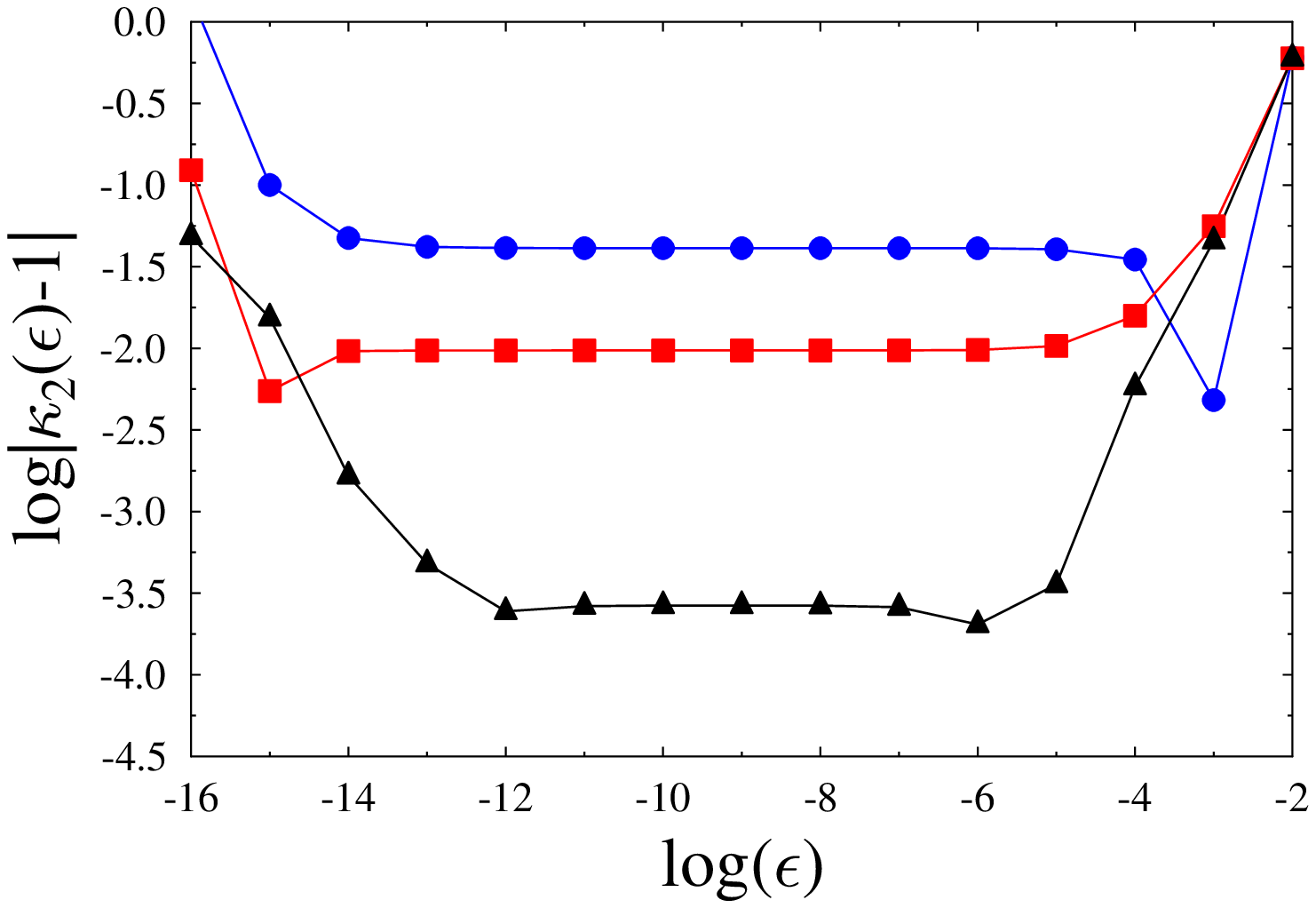}}}
\caption{Diagnostic quantities (a) $\log|\kappa_1(\epsilon)-1|$
  evaluated for Problem $P1$ and (b) $\log|\kappa_2(\epsilon)-1|$
  evaluated for Problem $P2$, cf.~\eqref{eq:kappa}, as functions of
  $\log\epsilon$ obtained using different discretizations of the time
  interval. {Blue circles: $N_T = 50$. Red squares: $N_T = 500$.
    Black triangles: $N_T = 5000$.} In all cases the perturbation
    direction is $g'_i = - r^3$, $i=1,2$.}
\label{fig:kappa}
\end{center}
\end{figure}

A key element of optimization algorithm \eqref{eq:desc} are the cost
functional gradients and a standard approach to their validation
consists in computing the directional G\^ateaux differential
$\J'_i(g_i;g'_i)$, $i=1,2$, for some arbitrary perturbations $g'_i$ in
two different ways, namely, using a finite-difference approximation
(with step size $\epsilon$) and using the inner product of the
adjoint-based gradient {with} the perturbation $g'_i$, namely Riesz
representation \eqref{eq:dJ}, and then examining the ratio of the two
quantities, i.e.,
\begin{equation}
\kappa_i (\epsilon) := \dfrac{\epsilon^{-1} \left[
\J_i(g_i + \epsilon g'_i) - \J_i(g_i) \right]}
{\int_{\I} \bnabla^{L_2} \J_i (r) \, g'_i(r) \, dr}, \quad i=1,2,
\label{eq:kappa}
\end{equation}
for a range of values of $\epsilon$. If the gradient $\bnabla^{L_2}
\J_i(r)$ is computed correctly, then for intermediate values of
$\epsilon$, $\kappa_i(\epsilon)$ will be close to the unity.
Remarkably, this behavior can be observed in Figures
\ref{fig:kappa}a,b corresponding to Problems $P1$ and $P2$ over a
range of $\epsilon$ spanning about 8 orders of magnitude. The quantity
shown in Figures \ref{fig:kappa}a,b is $\log| \kappa_i(\epsilon) - 1
|$, $i=1,2$, which represents the number of significant digits to
which the two ways to evaluate $\J'_i(g_i;g'_i)$ in \eqref{eq:kappa}
agree. Furthermore, we also observe that refining the resolution $N_T$
of the time interval $[0,T]$ yields values of $\kappa_i(\epsilon)$
closer to the unity. The reason is that in the
``optimize-then-discretize'' paradigm adopted here such refinement
of the discretization leads to a better approximation of the
continuous gradient \eqref{eq:gradJL2}. As can be expected, the
quantities $\kappa_i(\epsilon)$ deviate from the unity for very small
values of $\epsilon$, which is due to the subtractive cancellation
(round-off) errors in finite-differencing, and also for large values
of $\epsilon$, which is due to the truncation errors, both of which
are well-known effects.

\subsection{Computational Robustness of the Proposed Approach}
\label{ToC:Discussion_computational}
\begin{figure}
\begin{center}
\mbox{
\subfigure[]{\includegraphics[width=0.5\textwidth]{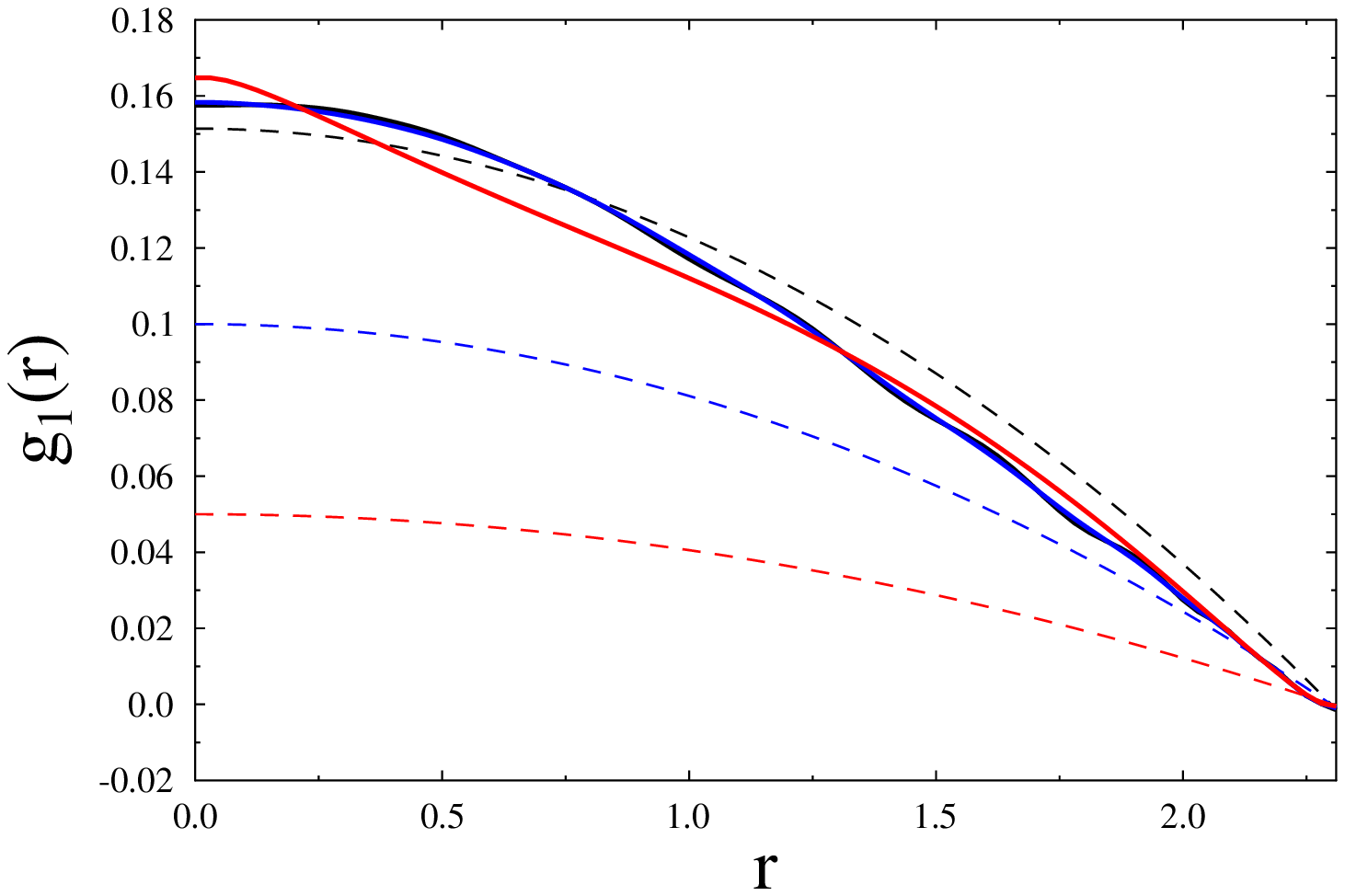}}\quad
\subfigure[]{\includegraphics[width=0.5\textwidth]{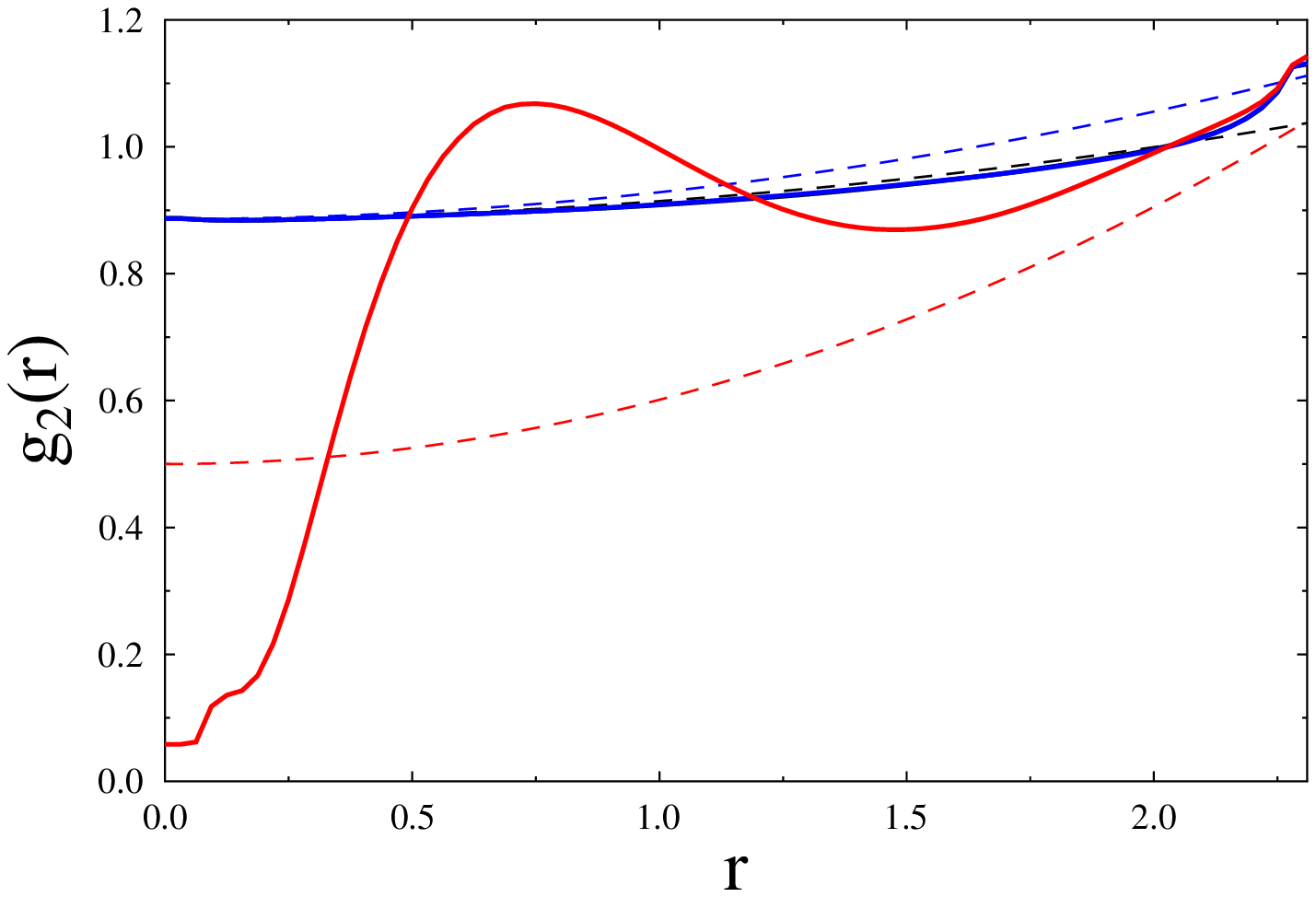}}}
\caption{{Reconstructions and the corresponding initial guesses
    in the solution of (a) Problem $P1$ and (b) Problem $P2$. Solid
    lines: reconstructions $\hat{g}_i$, $i=1,2$. Dashed lines: initial
    guesses $g_i^0$. Different} reconstructions are marked with the
  same color as the corresponding initial guesses.  Initial guesses
  and reconstructions discussed in Section \ref{ToC:Results} are
  marked in black.}
\label{fig:gig}
\mbox{
\subfigure[]{\includegraphics[width=0.5\textwidth]{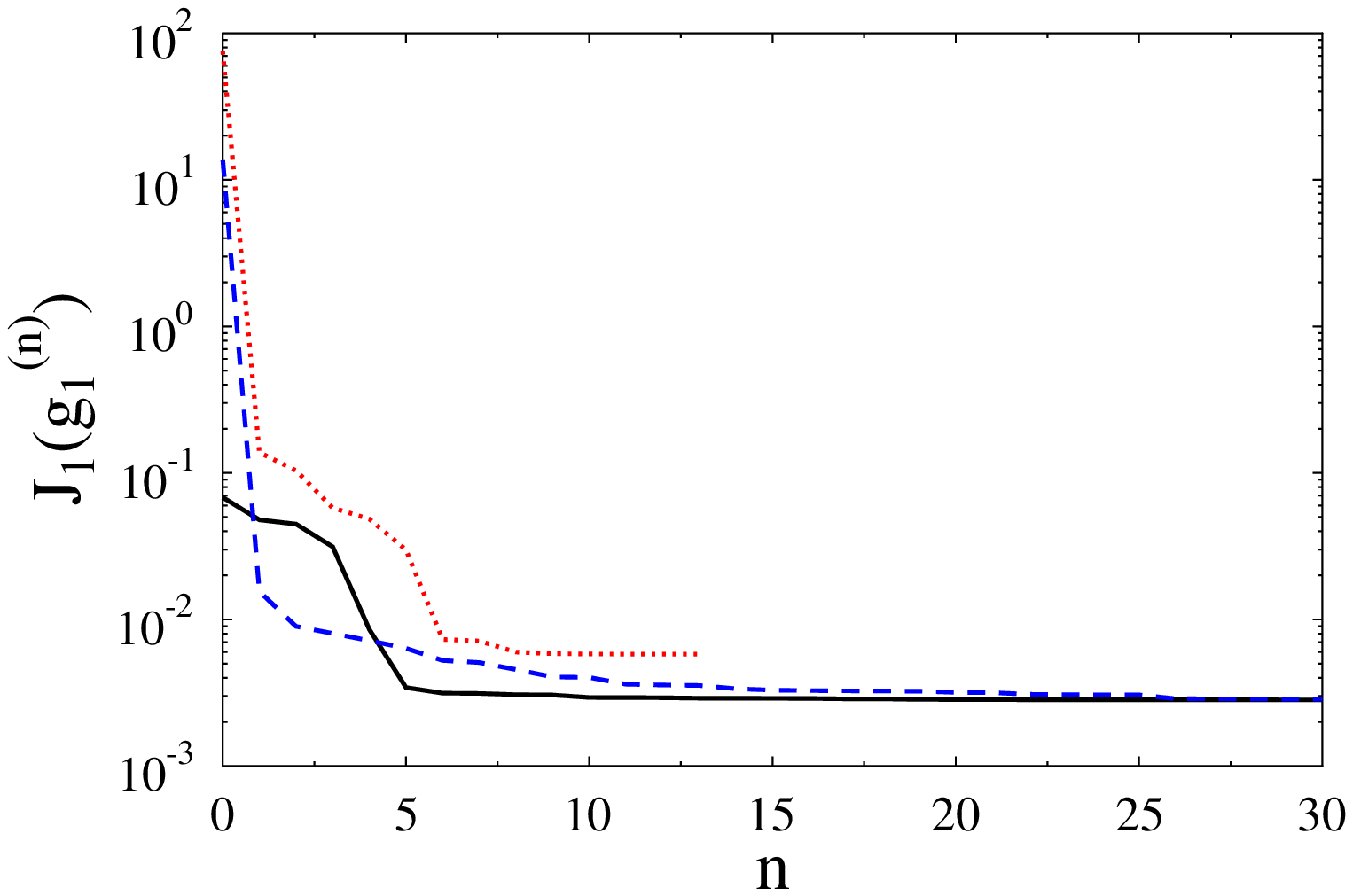}}\quad
\subfigure[]{\includegraphics[width=0.5\textwidth]{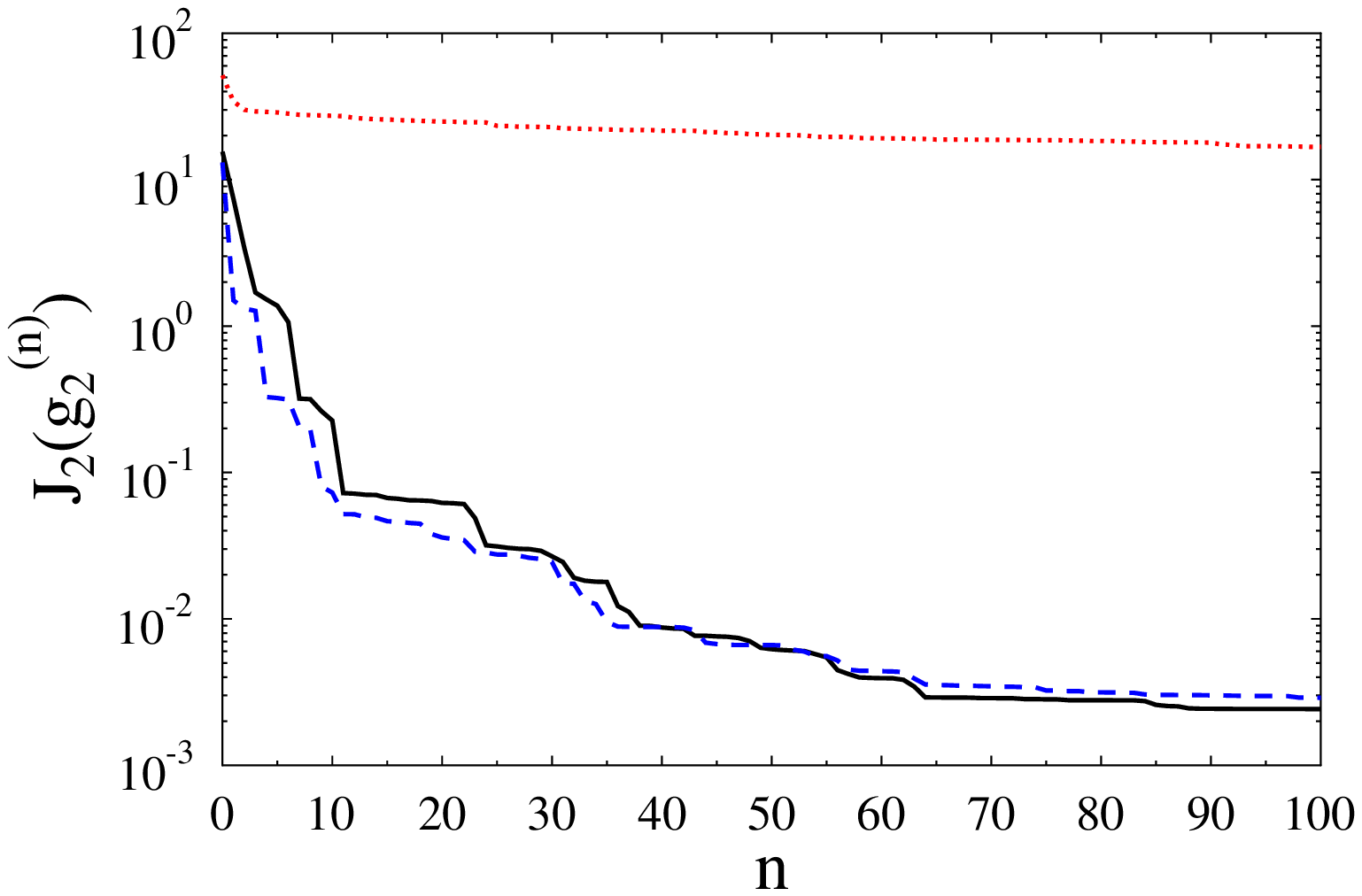}}}
\caption{Decrease of the cost functionals with iterations in (a)
  Problem $P1$ and (b) Problem $P2$ for the cases studied in Figure
  \ref{fig:gig} ({with the same color-coding}).
}
\label{fig:jig}
\end{center}
\end{figure}

In this Section we focus on the effect that the choice of initial
guess $g_i^0$, $i=1,2$, has on the reconstructed functions $\hat{g}_1$
and $\hat{g}_2$. We note that, given the nonlinearity of system
\eqref{eq:DescriptorSystemB}, optimization problems $P1$ and $P2$ may
be nonconvex and optimality conditions \eqref{eq:opt} and
\eqref{eq:dJ2} characterize minimizers which are only {\em local}.
Thus, different initial guesses may in principle give rise to
different reconstructions and convergence to a global minimum cannot
be a priori assured. We investigate this issue in Figures
\ref{fig:gig}a and \ref{fig:gig}b where we show the reconstructions
obtained, respectively, in Problems $P1$ and $P2$ using {different}
initial guesses generally much worse than mean-field model
\eqref{eq:IG} used in Section \ref{ToC:Results}. As regards Problem
$P1$, we note in Figure \ref{fig:gig}a that accurate reconstructions
are obtained using even relatively poor initial guesses $g_1^0$. On
the other hand, in Figure \ref{fig:gig}b we see that in Problem $P2$
the reconstruction fails for a less accurate initial guess $g_2^0$.
These two examples are representative of the behavior we generally
observed in our calculations and we conclude that Problem $P1$ appears
more robust with respect to the choice of the initial guess than
Problem $P2$. We also noted that in both problems convergence tends to
be more sensitive to the values assumed by the initial guesses $g_1^0$
and $g_2^0$ at $r=0$ and $r=r^\circ$ than to their behavior for
intermediate values of $r$.  Results from Figure \ref{fig:gig} are
corroborated by the corresponding histories of the cost functionals in
Figures \ref{fig:jig}a and \ref{fig:jig}b.  We note that in the case
of the poorest initial guess in Problem $P2$, cost functional
$\J_2(g_2^{(n)})$ reveals hardly any decrease with the iterations at
all. While in the cases of successful reconstructions the cost
functionals $\J_1(g_1^{(n)})$ and $\J_2(g_2^{(n)})$ drop over several
orders of magnitude, they never attain values lower than
$\O(10^{-3})$. This is a consequence of the phase-dependent behavior
of the measurements which cannot be resolved using ansatz in the form
\eqref{eq:DescriptorSystemB}, cf.~Assumption \ref{ass1}(a), see also
the inset in Figure \ref{fig:ra1a3}a.

Inverse problems of the type considered here often tend to be
ill-posed, in the sense that small perturbations to the data, for
example due to noise, may result in significant changes in the
computed solution. Suitable regularization, for instance, using
Tikhonov's technique \cite{t05,v02}, may be required to stabilize the
solution procedure in such situations. To focus attention in the
present study, we concentrated on the structure of the gradients and
did not investigate the effect of noise on the reconstructions, hence
such regularization was not necessary. We refer the reader to
\cite{bvp10,bp11a} for a thorough analysis of regularization applied to a
related reconstruction problem.

%***********************************************************************
\subsection{Physical Interpretation of the Results}
\label{ToC:Discussion_physical}

In this Section we {propose} some physical interpretation of the
numerical reconstruction results from Section \ref{ToC:Results}.  We
note in Figures \ref{fig:g123} and \ref{fig:ra1a3} that the identified
phase-invariant oscillation model \eqref{eq:DescriptorSystemB} is in
fact in remarkably good agreement with the data obtained from the
solution of the Navier-Stokes equation \eqref{Eqn:IncompressibleFlow}.
A small difference between the model and the data is visible as
wiggles due to the second harmonic {present in the measurements} which
violates the phase-invariance assumed in our model ansatz,
cf.~Assumption \ref{ass1}(a).  This difference can be easily removed
by a simple pre-processing of the measurement data.  Referring to
Appendix \ref{ToC:POD}, we note that the POD eigenvalue $\lambda_1$,
representing the variance of $a_1$, is larger than the {eigenvalue
  $\lambda_2$ which represents the variance of $a_2$}.  The following
rescaling transformation
\begin{subequations}
\label{eqn:a1a2rescaling}
\begin{eqnarray}
{\bar{a}_1} &=& \sqrt{\frac{\lambda_1+\lambda_2}{2 \lambda_1}} \> a_1,
\\
{\bar{a}_2} &=& \sqrt{\frac{\lambda_1+\lambda_2}{2 \lambda_2}} \> a_2
\end{eqnarray}
\end{subequations}
ensures equipartition of energy in the new variables $\bar{a}_1$ and
$\bar{a}_2$ while conserving the total energy in both modes.  This
rescaling effectively removes the second harmonics from $\bar{a}_i$,
$i=1,2$, which could be used as new inputs {for the reconstruction}.
However, we refrained from applying \eqref{eqn:a1a2rescaling} in the
computations reported in Section \ref{ToC:Results} in order to show
the power of the proposed identification method to deal with data
which cannot be perfectly matched by the model.

Another observation concerning Figure \ref{fig:g123} is the
significant deviation of reconstructed functions $\hat{g}_1$ and
$\hat{g}_2$ from the parabolic mean-field relations
\eqref{Eqn:MeanFieldSystem}.  Evidently, higher-order corrections,
such as {$r^4$, $r^6$, etc.}, are required for a better agreement
between the identified propagators $\hat{g}_i$, $i=1,2$, and the
expansions used in the mean-field model. The information which
higher-order terms ought to be included in the model as opposed to an
a priori fixed polynomial expansion used typically in model
identification is therefore the unique advantage of the proposed
identification strategy.  We note that odd powers of $r$ can be
excluded by phase-invariance considerations. In addition, with our
reconstruction method we were able to determine more accurate values
of the instability growth rate at the origin and the oscillation
frequency at the limit cycle than used in mean-field model
\eqref{eq:IG}. We stress that in fact such seemingly insignificant
modifications of the structure of the reduced-order model may already
affect its utility for various control applications.

These results also shed light on the validity of mean-field model
\eqref{Eqn:MeanFieldSystem}.  Initially, the mean-field theory
\cite{Stuart1958jfm,Watson1960jfm} was derived to be valid near the
onset of the oscillation only.  We probed the applicability of this
model by applying it at a Reynolds number $100$ which is more than
twice the critical value of $47$.  Hence, the deviation of $\hat{g}_1$
and $\hat{g}_2$ from \eqref{Eqn:LandauEquation} does not invalidate
the mean-field theory.  One reason for this deviation is the change
of the structure of the vortex street during the transient.  The
optimal oscillatory modes deform from {the stability} eigenmodes
{into the} POD modes while the fluctuation center moves upstream
and the frequency and wavenumber increase
\cite{Noack2003jfm,Tadmor2011ptrsa}.  Similarly, the mean-field
correction (the shift mode,
{cf.~\eqref{Eqn:MeanFieldExpansion_2}}) changes during the
transient \cite{Tadmor2010pf} {which} has a noticeable effect on
the  mean-field model
\cite{Morzynski2006aiaa}.

Finally, the results concerning the identified descriptor system are
{also} relevant to the empirical 9-dimensional Galerkin model
accounting for the base-flow variation and for the first four
harmonics \cite{Noack2003jfm}.  The initial exponential growth of the
first harmonic is limited by the base-flow variation (which reduces
the production of fluctuation energy) and, to a lesser extent, by the
energy transfer from the first into higher harmonics.  The energy
transfer may be accounted for by an energy-dependent eddy viscosity
in the mean-field system.  Under certain assumptions (see, e.g.\
\cite[chap.\ 3]{Noack2011book}), a generalized Landau equation
\begin{equation}
\dot r = \sigma_1 r - \beta r^3 - \gamma r^3
\label{eq:Landau9D}
\end{equation}
can be derived, where $\sigma_1$ denotes the growth rate near the
fixed point $r=0$ and $\beta$, $\gamma$ characterize the damping from
the $0$-th and from higher harmonics, respectively. By carefully
comparing the 9-dimensional Galerkin model with the identified
phase-invariant system {it} may be therefore {possible to determine
  the values of} $\beta$ and $\gamma$, or even to correct the powers
of the new terms in \eqref{eq:Landau9D}.  A complete derivation and an
in-depth discussion of this problem is outside the scope of the
present study.

%***********************************************************************
\section{Conclusions and Future Directions}
\label{ToC:Conclusions}

We have proposed and validated a novel method for model identification
which is an adaptation of an approach already used in the context of
systems described by PDEs \cite{bvp10,bp11a}. {As indicated in
  Figure \ref{fig:f},} we depart from the traditional approach of (1)
characterizing the propagator of the dynamical system in a parameter
space and {then} (2) performing a parameter identification.  Thus,
arbitrary polynomial expansions of the propagator may be performed a
posteriori ({following the solution of the optimal reconstruction
  problem}) at a practically vanishing cost.  In addition, the
performance of parametric models may easily be assessed {and}, if
necessary, improved by {introducing} additional terms {motivated by
  the form of the reconstructed constitutive relation}.

The method is applied to a three-dimensional descriptor system with
three {a priori undetermined relations} describing the fluctuation
growth, frequency and mean-field correction as functions of the
fluctuation energy.  As a benchmark problem {we chose} the onset of
laminar von K\'arm\'an vortex shedding behind a circular cylinder.
Results of a direct numerical simulation are transcribed into the mode
amplitudes of a minimal 3-state Galerkin model \cite{Noack2003jfm}
which are then captured with remarkable accuracy by our identified
descriptor system. The form of the reconstructed system is marked by a
noticeable departure from the mean-field models and may therefore
guide the refinement of the latter by inclusion of higher-order terms.
We emphasize that the usefulness of reduced-order models for flow
control applications may in fact depend on such differences.

As regards future research directions, while the present results offer
a proof of the concept for the proposed approach based on a rather
well-understood example, the key question is extension of this method
to the identification of models with more complicated structure
featuring, for example, {multiple time scales}, state space of a
higher dimension, coexistence of several oscillation frequencies,
non-trivial phase dependence and higher-dimensional inertial
manifolds. {As regards the first issue, one can consider a
  modification of our model problem
  \eqref{eq:DescriptorSystemB}--\eqref{eq:DescriptorSystemA} with
  Assumption \ref{ass1}(b) revised to allow $a_3(t)$ to be a ``fast''
  variable. While in such setting our computational approach would
  formally remain unchanged (except that Problem $P3$ would be
  replaced with a problem similar to $P1$ or $P2$), it is interesting
  how it would actually perform in practice.}  Dealing with some of
{the other} aspects will require formulation of the
reconstruction problems in terms of propagator functions depending on
more than just one state variable ($r$ in the examples considered in
the present study). This will, in turn, lead to a number of
interesting questions at the level of numerical analysis and
scientific computing related to the evaluation of the cost functional
gradients. {An emerging application which involves some of the
  aforementioned extensions is related to the question of optimal
  parametrization of subgrid turbulence representations which is an
  important open problem in theoretical fluid mechanics \cite{lm99}.
  In the context of Galerkin reduced-order models, it may take the
  form of an additional dissipative term with the magnitude
  proportional to an ``eddy viscosity'' $\nu_T$
\begin{equation}
\frac{d\bm{a}}{dt} = \bm{f}(\bm{a}) + \nu_T \, \bL^{\nu}\, \bm{a}, \quad \bm{a} \in \RR^N,
\label{eq:nuT}
\end{equation}
}{where $\nu_T \, \bL^{\nu}\, \bm{a}$ represents the stabilizing
  viscous term of the Navier-Stokes equation in the Galerkin system.
  In general, this term can be proven to be energy dissipative for all
  orthonormal systems of modes and a large class of boundary
  conditions. For some analytical modes, e.g., the Stokes modes, it
  can be shown that matrix $\bL^{\nu}$ is diagonal and
  negative-definite.}  {There is abundant evidence
  \cite{Noack2011book,wabi12,Balajewicz2013jfm} that nonlinear closure
  strategies perform better as regards stabilization of system
  \eqref{eq:nuT}.  Assuming $\nu_T = \nu_T(\| \bm{a} \|)$ gives rise
  to an identification problem analogous to $P1$ and $P2$, and one can
  use the algorithm described in Section \ref{ToC:Computation} to
  determine optimal closure strategies leading to the best possible
  reconstruction of the available data. Preliminary identification
  results already obtained based on a reduced-order model
  \eqref{eq:nuT} with dimension $N=20$ applied to a complex
  mixing-layer flow are quite encouraging and reveal some nontrivial
  physical insights. They will be reported in the near future upon
  completion of the study.}

We also remark that a surprisingly large set of modelling and control
problems can be cast in a similar form of function identification of a
descriptor system
\begin{equation}
\frac{d\bm{a}}{dt} = \bm{f}(\bm{a},\bm{b}), \quad
\bm{b} = \bm{g} (\bm{a}) .
\end{equation}
For reasons of simplicity, let us assume that function $\bm{f}$ is
known and that function $\bm{g}$ needs to be determined.  If $\bm{b}$
characterizes the slow modes, then $\bm{b} = \bm{g}( \bm{a})$
represents the inertial manifold to be {identified} from a given
{system} trajectory \cite{Gorban2005book}.  If $\bm{b}$ represents
high-frequency components, such as the parameters of a subgrid
turbulence representation {described above}, then their
functional dependence on the state {variable} $\bm{a}$ may also be
inferred with our approach.  {On the other hand,} if $\bm{b}$ denotes
the actuation amplitudes, as in numerous wake flow stabilization
studies
\cite{Gerhard2003aiaa,Bergmann2005pf,Thiria2006jfm,Pastoor2008jfm,Weller2009physd},
then $\bm{b}=\bm{g}(\bm{a})$ represents a full-state {feedback}
control law.  {In principle,} this control law may {as well} be
identified from desired trajectories $t \mapsto \bm{a}(t)$.

\section*{Acknowledgements}
The authors acknowledge the funding 
and excellent working conditions 
of the  Chair of Excellence
'Closed-loop control of turbulent shear flows 
using reduced-order models' (TUCOROM)
of the French Agence Nationale de la Recherche (ANR)
and hosted by Institute PPRIME.
The first author is, in particular, grateful for the hospitality 
of this Chair of Excellence at Institute PPRIME where most of this work was carried out.
We also thank the Ambrosys Ltd.\ Society for Complex Systems Management
and the Bernd Noack Cybernetics Foundation for additional support.
We appreciate valuable stimulating discussions
with Markus Abel, Robert Niven, Michael Schlegel and Gilead Tadmor
as well as the local TUCOROM team:
Jean-Paul Bonnet,
Laurent Cordier,
Thomas Duriez,
Peter Jordan, 
Vladimir Parezanovic and
Andreas Spohn.
%Last but not least, we thank the referees
%for their thoughtful and helpful suggestions.
%\end{acknowledgments}
%***********************************************************************

%\clearpage

%----- Appendix --------------------------------------------------------
\begin{appendix}
\section{Proper Orthogonal Decomposition}
\label{ToC:POD}

In this Appendix we describe the Proper Orthogonal Decomposition (POD)
employed in Section \ref{ToC:GalerkinExpansion} to construct a
low--dimensional model from the simulation data.  {It is closely
  related to other techniques of data analysis known as the Principal
  Component Analysis or, in the discrete setting, the Singular--Value
  Decomposition.}  The starting point are snapshots of the velocity
field $\bm{u}^m(\bm{x}) {:=\u(\x,t^m)}$, $m=1,\ldots,M$, $\bm{x}
\in \Omega$, cf.~\eqref{eq:Omega}.  These snapshots are sampled
{at times $t^m$ uniformly spaced over} one period of oscillation
and form a statistically representative ensemble for the considered
first and second moments.

The goal is to construct a 'least--order' Galerkin expansion 
%----- Equation --------------------------------------------------------
\begin{equation}
\label{eqn:Galerkin_expansion}
\bm{u}(\bm{x},t) = \bm{u}_0 (\bm{x}) + \sum\limits_{i=1}^N a_i(t) \> \bm{u}_i(\bm{x}) + \bm{u}_{\rm res}(\bm{x},t)
\end{equation}
%-----------------------------------------------------------------------
with base mode $\bm{u}_0$ {and} $N$ space-dependent expansion
modes (basis functions) $\bm{u}_i(\bm{x})$ {with} the
corresponding mode amplitudes $a_i(t)$ which will result in a
minimum-norm average residual $\bm{u}_{\rm res}$ of the snapshot
ensemble (see, e.g., \cite{Fletcher1984book, Noack2011book,
  Holmes2012book}).  The base mode is necessary so that Galerkin
expansion \eqref{eqn:Galerkin_expansion} satisfies inhomogeneous
boundary conditions for arbitrary mode amplitudes. {For instance,
  expansion \eqref{eqn:Galerkin_expansion} captures the prescribed
  oncoming flow velocity regardless of the values of $a_i$.}

The Galerkin expansion and its residual are embedded in the Hilbert
space $L_2( \Omega )$ of square-integrable vector fields.  The inner
product of two elements $\bm{v}, \bm{w} \in L_2( \Omega )$ is defined
by
%----- Equation --------------------------------------------------------
\begin{equation}
{ \big\langle \bm{v}, \bm{w} \big \rangle_{{L_2(\Omega)}} := \int\limits_{\Omega} \, \bm{v} \cdot \bm{w} \, d\bm{x}},
\end{equation}
%-----------------------------------------------------------------------
where '$\cdot$' denotes the standard Euclidean inner product and
$d\bm{x}$ an infinitesimal volume element of the domain $\Omega$.  The
associated norm thus is
%----- Equation --------------------------------------------------------
\begin{equation}
\Vert \bm{u} \Vert_{{L_2(\Omega)}} := \sqrt{ \langle \bm{u}, \bm{u} \rangle_{{L_2(\Omega)}} }.
\end{equation}
%-----------------------------------------------------------------------
We search for empirical modes $\bm{u}_i$, $i=0,\ldots,N$, where $N \le
M-1$, which minimize the average residual of the Galerkin expansion of
the snapshots
%----- Equation --------------------------------------------------------
\begin{equation}
\label{eqn:GE_Snapshots}
\bm{u}^m  = \bm{u}_0 + \sum\limits_{i=1}^N a_i^m \> \bm{u}_i + \bm{u}_{\rm res}^m,
\quad m=1,\ldots,M,
\end{equation}
%-----------------------------------------------------------------------
with optimal mode amplitudes $a_i^m$ in the sense of the $L_2$ norm, i.e.,
%----- Equation --------------------------------------------------------
\begin{equation}
\label{eqn:MinimumResidual}
\overline{ \left\Vert \bm{u}_{\rm res} \right\Vert^2_{{L_2(\Omega)}} }
:= \frac{1}{M} \sum\limits_{m=1}^M \left \Vert \bm{u}_{\rm res}^{m} \right\Vert^2_{{L_2(\Omega)}}
{=} \hbox{min}.
\end{equation}
%-----------------------------------------------------------------------
This problem is solved by the snapshot POD \cite{Sirovich1987qam1}
{and} the base mode is the mean of the snapshots
\begin{equation}
\bm{u}_0 := \frac{1}{M} \sum\limits_{m=1}^M \bm{u}^m.
\end{equation}
The POD modes arise from the correlation matrix
$\bm{C} := \left( C^{mn} \right)_{m,n=1,\ldots,M}$
of the snapshot fluctuations
\begin{equation}
C^{mn} := \frac{1}{M} \Big\langle\bm{u}^m - \bm{u}_0, 
                            \bm{u}^n - \bm{u}_0 \Big\rangle_{{L_2(\Omega)}}.
\end{equation}
We note that $\bm{C}$ is a positive semi-definite Grammian matrix.
Hence, the eigenvalue problem
\begin{equation}
\bm{C} \bm{e}_i = \lambda_i \bm{e}_i, \quad i=1,\ldots,M,
\end{equation}
yields an orthonormal set of real eigenvectors $\bm{e}_i =
\left[e_i^1,\ldots,e_i^M \right]^T $, $\bm{e}_i^T
\bm{e}_j = \delta_{ij}$, with non-negative eigenvalues which can be
ordered as
\begin{equation}
\lambda_1 \ge \lambda_2 \ge \ldots \ge \lambda_M = 0.
\end{equation}
The last equality arises from the fact that $M$ vectors span a
subspace of maximum dimension $M-1$.  Hence, $M$ snapshots define only
$M-1$ POD modes and the corresponding amplitudes. These are given by
\begin{equation}
\label{eqn:POD}
\bm{u}_i = \frac{1}{\sqrt{M \lambda_i}} 
\sum\limits_{m=1}^M e_i^m \> \left( \bm{u}^m - \bm{u}_0 \right),\quad
a_i^m = \sqrt{\lambda_i M} \> e_i^m,
\quad i=1,\ldots,M-1.
\end{equation}
The POD modes form an orthonormal basis in $L_2(\Omega)$, $\langle
  \bm{u}_i, \bm{u}_j \rangle_{{L_2(\Omega)}} = \delta_{ij}$,
$i,j=1,\ldots,N$, while their amplitudes have vanishing means
$\overline{a_i} = 0$, $i=1,\ldots,N$, and diagonal second moments
$\overline{a_i \> a_j} = \lambda_i \> \delta_{ij}$, $i,j=1,\ldots,N$.
The average is to be understood in terms of the snapshot ensemble
(see, e.g., \eqref{eqn:MinimumResidual}).  The eigenvalue $\lambda_i$
can be interpreted as twice the fluctuation energy contained the
$i$-th mode.
%***********************************************************************
\section{Regularity of Reconstructed Function $g_1$ versus Existence
  and Uniqueness of Solutions to Equation \eqref{eq:DescriptorSystemB_1}}
\label{ToC:Regularity}

By the one-dimensional embedding result $H^1(\I) \rightharpoonup
C^{0,\frac{1}{2}}(\I)$ \cite{af05}, we note that the reconstructed
function will be H\"older-continuous with $\lambda=1/2$. Thus, it
will not meet the assumptions of the Picard-Lindel\"of theorem
\cite{mm82}, and in principle will ensure existence only, without
uniqueness, of solutions of equation \eqref{eq:DescriptorSystemB_1}.
In order to ensure the Lipschitz-continuity ($\lambda=1$) of $g_1$,
which would also guarantee uniqueness of solutions of
\eqref{eq:DescriptorSystemB_1}, we would need to reconstruct $g_1$ as
an element of Sobolev space $H^2(\I)$, because then $H^2(\I)
\rightharpoonup C^{0,1}(\I)$. While there are no fundamental
difficulties here (we would need to replace inner product
\eqref{eq:ipH1} with the corresponding definition in $H^2(\I)$), we
will refrain from this in the actual computations in Section
\ref{ToC:Results} in order to keep the approach as simple as possible.
Nevertheless, in all problems we treated with the proposed approach
the reconstructed functions possessed the Lipschitz regularity which
was verified a posteriori by performing suitable grid-refinement
studies.

%***********************************************************************
\section{Derivation of Perturbation Equation \eqref{eq:dxi1}}
\label{ToC:Perturbation}

In this Appendix we present a derivation of perturbation equation
\eqref{eq:dxi1}. Obtaining this equation is made somewhat more
involved by the fact that the perturbation variable $g_1$ is itself a
function of the state magnitude, i.e., $g_1 = g_1(r)$. We assume here
that only $g_1$ is perturbed while $g_2$ remains fixed, with the
opposite case leading to essentially the same calculations. By
substituting, respectively, $g_1 = g_{1a}$ and $g_1 = g_{1b}$ into
equation \eqref{eq:DescriptorSystemA_1}, we obtain $\dot{\bxi}_a =
\left( g_{1a}(r_a) \, \bI + g_{2}(r_a) \, \bJ \right)\, \bxi_a$ and
$\dot{\bxi}_b = \left( g_{1b}(r_b) \, \bI + g_{2}(r_b) \, \bJ
\right)\, \bxi_b$ , where $\bxi_a := \bxi(g_{1a})$ and $\bxi_b :=
\bxi(g_{1b})$ are the corresponding solutions and $r_a := | \bxi_a |$,
$r_b := | \bxi_b |$. Taking the difference of these two equations and
defining $\bxi' := \bxi_a - \bxi_b$. we obtain
\begin{equation}
\begin{aligned}
\dot{\bxi}' &=   g_{1a}(r_a) \, \bI\,\bxi_a + g_{2}(r_a) \, \bJ \, \bxi_a - g_{1b}(r_b) \, \bI\,\bxi_b - g_{2}(r_b) \, \bJ \, \bxi_b \\
 &=   g_{1a}(r_a) \, \bI\,\bxi_a + g_{2}(r_a) \, \bJ \, \bxi_a + g_{1b}(r_a) \, \bI\,\bxi_a - g_{1b}(r_a) \, \bI\,\bxi_a
- g_{1b}(r_b) \, \bI\,\bxi_b - g_{2}(r_b) \, \bJ \, \bxi_b \\
& = g'(r_a) \, \bI \, \bxi_a + \underbrace{g_{1b}(r_a) \, \bI\,\bxi_a - g_{1b}(r_b) \, \bI\,\bxi_b}_A +
\underbrace{g_{2}(r_a) \, \bJ\,\bxi_a - g_{2}(r_b) \, \bJ\,\bxi_b}_B,
\end{aligned}
\label{eq:dxi3}
\end{equation}
where we also set $g'(\cdot) = g_{1a}(\cdot) - g_{1b}(\cdot)$ in the
first term on the RHS. As regards the terms denoted $A$ in
\eqref{eq:dxi3}, they are transformed as follows using the fundamental
theorem of calculus for line integrals and the change of variables
$\bxi(s) = \bxi_b + s \, (\bxi_a - \bxi_b)$ for $s \in [0,1]$
\begin{equation}
\begin{aligned}
A = g_{1b}(r_a) \, \bI\,\bxi_a - g_{1b}(r_b) \, \bI\,\bxi_b 
& = \int_{\bxi_b}^{\bxi_a} \bnabla_{\bxi} (g_{1b}(\bxi) \, \bxi)\, d\bxi \\
& = \left( \int_0^1 \bnabla_{\bxi} \bF(\bxi_b + s \bxi') \, ds \right) \, \bxi',
\end{aligned}
\label{eq:A}
\end{equation}
where we also denoted $\bF(\bxi) := g_{1b}(\bxi) \, \bxi$. The
integrand expression on the RHS in \eqref{eq:A} is then expanded in
the Taylor series around $\bxi' = \0$
\begin{equation}
\Dpartial{}{\xi_j} F_i(\bxi_b + s\bxi') = \Dpartial{}{\xi_j}F_i(\bxi_b) + \Dpartialmix{}{\xi_k}{\xi_j}F_i(\bxi_b)\, \xi'_k \, s + \O\left( |\bxi'|^2 \right),
\quad i,j,k=1,2,
\label{eq:Ftay}
\end{equation}
where the component notation was used for clarity with $F_i$, $\xi_j$
and $\xi'_k$ denoting the components of vectors $\bF$, $\bxi$ and
$\bxi'$. Plugging expansion \eqref{eq:Ftay} into expression
\eqref{eq:A} we obtain
\begin{equation}
\begin{aligned}
A & = \bnabla_{\bxi} (g_{1b}(\bxi) \, \bxi )\Big|_{\bxi=\bxi_b}\, \bxi' \, \int_0^1 \, ds + \O\left( |\bxi'|^2 \right) \\
& = \left[ g_{1b}(r_b) \, \bI + \bI \, \bxi_b \,\left(\bnabla g_{1b}(r_b)\right)^T \right] \, \bxi' + \O\left( |\bxi'|^2 \right).
\end{aligned}
\label{eq:A2}
\end{equation}
Noting that term $B$ in \eqref{eq:dxi3} transforms in an analogous way
to $A$ in \eqref{eq:A}--\eqref{eq:A2}, using these results in
\eqref{eq:dxi3}, assuming smallness of $\bxi'$ and dropping terms of
order quadratic and higher we finally arrive at perturbation equation
\eqref{eq:dxi1}.
\end{appendix} 
%----- Main part -------------------------------------------------------

%---- Bibliography -----------------------------------------------------
\bibliographystyle{unsrt}
%\bibliography{~/texts/ROMID/Protas_pre/Bib/all}
%\bibliography{/Users/bprotas/texts/ROMID/Protas_pre/Bib/all}
%\bibliography{/home/bprotas/texts/ROMID/Protas_pre/Bib/all}
%\bibliography{./allPROTAS,./allNOACK}

\end{document}